\documentclass[a4paper, 12pt]{article}

\usepackage[ansinew]{inputenc}
\usepackage[english]{babel}
\usepackage{textcomp}
\usepackage{csquotes}
\usepackage{mathrsfs}
\usepackage{amsfonts}
\usepackage{mathtools}
\usepackage{amssymb}
\usepackage{amsmath}
\usepackage{amsthm}
\usepackage{amscd}
\usepackage{tensor}
\usepackage{stmaryrd}
\usepackage{rotating}
\usepackage[text={6.8in,8.5in},centering]{geometry}
\usepackage{tikz}
\usetikzlibrary{knots}
\usetikzlibrary{arrows, calc, decorations.markings, positioning}
\usepackage{tikz-cd}
\usepackage{relsize} 
\usepackage[none]{hyphenat}
\usepackage{epigraph}
\usepackage{tikz-cd}

\tikzcdset{row sep/normal=0.6cm}
\tikzcdset{column sep/normal=0.2cm}

\tikzcdset{row sep/small=0.3cm}
\tikzcdset{column sep/large=0.6cm}

\tikzcdset{row sep/large=0.9cm}
\tikzcdset{column sep/small=0.1cm}

\usepackage{graphicx}
\graphicspath{ {./Images/} }

\newtheorem{definition}{Definition}[section]

\newcommand{\Real}{\mathbb{R}}

\usepackage{hyperref,xcolor}
\usepackage{makeidx}
\definecolor{nred}{rgb}{0.7,0,0}
\definecolor{nblue}{rgb}{0,0,0.6}
\definecolor{ngreen}{rgb}{0,0.6,0}

\hypersetup{
    colorlinks=true,
    linkcolor=nred,
    filecolor=nred,      
    urlcolor=ngreen,
    citecolor=nblue,
}

\usepackage[
backend=biber,
style=alphabetic,
sorting=ynt
]{biblatex}
\begin{document}

\title{\textbf{An Invitation to Higher Arity Science}}

\author{Carlos Zapata-Carratal\'a$^{1, }$\footnote{lead author, \url{c.zapata.carratala@gmail.com}}  {\,\,\,}   Xerxes D. Arsiwalla$^{2, 3, }$\footnote{\url{x.d.arsiwalla@gmail.com}}   \\
{}  \\
{\it \small $^{1}$School of Mathematics, The University of Edinburgh, United Kingdom}\\
{\it \small $^{2}$Pompeu Fabra University, Barcelona, Spain}\\
{\it \small $^{3}$Wolfram Research, USA}    
}

\date{}

\maketitle

\sloppy

\begin{abstract}
Analytical thinking is dominated by binary ideas. From pair-wise interactions, to algebraic operations, to compositions of processes, to network models, binary structures are deeply ingrained in the fabric of most current scientific paradigms. In this article we introduce arity as the generic conceptualization of the order of an interaction between a discrete collection of entities and argue that there is a rich universe of higher arity ideas beyond binarity waiting to be explored. To illustrate this we discuss several higher order phenomena appearing in a wide range of research areas, paying special attention to instances of ternary interactions. From the point of view of formal sciences and mathematics, higher arity thinking opens up new paradigms of algebra, symbolic calculus and logic. In particular, we delve into the special case of ternary structures, as that itself reveals ample surprises: new notions of associativity (or lack thereof) in ternary operations of cubic matrices, ternary isomorphisms and ternary relations, the integration problem of 3-Lie algebras, and generalizations of adjacency in 3-uniform hypergraphs. All these are open problems that strongly suggest the need to develop new ternary mathematics. Finally, we comment on potential future research directions and remark on the transdisciplinary nature of higher arity science.
\end{abstract}

\textbf{keywords} - higher order structures, complex systems, interdisciplinary research, generalized associativity, ternary algebra, hypergraphs, matrix algebra, category theory

\clearpage

\tableofcontents

\newpage

\epigraph{The task is not to see what has never been seen before, but to think what has never been thought before about what we see everyday.}{Erwin Schr\"odinger}

\section{Introduction} \label{intro}

Discrete bits of quantifiable and communicable information -- words in speech and text, digits in computers and records, and symbols in formulae and algorithms -- are the building blocks of the human intellectual edifice.  However, the architectural integrity of this construction relies on much more than the mere collection of its parts: basic ideas are connected, organized and composed thanks to the cohesive power of conceptual relations and associations. Our aim is to bring this structural aspect of modern science to the foreground and to illustrate how our understanding of its underpinnings may lead to the development of new theories and paradigms.\newline

To this end we shall introduce the general concept of `arity', which refers to the order or degree of a relation between a family of countable entities. For example, a predator-prey dynamic in an ecosystem is an instance of arity 2, i.e. a \emph{binary} relation between two species, and the chemical bond of a methane molecule is an instance of arity 5, i.e. a \emph{quinary} relation between five atoms. Motivating our discussion with elementary examples of ternary structures in mathematics, we will argue that restricting ourselves to the use of binary ideas, which currently dominate the intellectual landscape, has a limiting effect on our capacity to develop new conceptual paradigms. We will also present several ternary (and higher) phenomena appearing in a wide range of disciplines to illustrate the untapped potential of higher arity science\footnote{A short video summary by one of the authors can be found here: \url{https://www.youtube.com/watch?v=62UFbGsj5Jg}.}. For reasons we shall discuss below, the full scope of arity, particularly in theoretical frameworks involving higher order structures, is yet to be realized.\newline

In recent years, an increasing number of voices coming from multiple research communities, including network science \cite{courtney2016generalized,battiston2020networks}, mathematical physics \cite{kerner2008ternary,azcarraga2010nary,baez2011invitation,benini2021operads}, biology \cite{klamt2009hypergraphs,kempes2021multiple}, neuroscience \cite{zhou2006learning,yu2011higher}, ecology \cite{billick1994higher,mayfield2017higher}, complexity science \cite{baas2009new,courtney2016generalized,neuhauser2021consensus} and computer science \cite{zhou2006learning,moskovich2015tangle,wolfram2021multicomputation}, are echoing concerns about the shortcomings of mainstream theoretical and computational paradigms due to their limited binary nature. This trend was succinctly captured in a recent Quanta Magazine article \cite{quanta2021higher}. The preponderance of binarity, and its potential limiting consequences to the development of new science, can only be adequately understood in the context of a general framework integrating higher forms of arity. This article aims to provide a first step towards such a framework.\newline

From a purely mathematical point of view, higher arity has been largely historically ignored -- almost to a surprising degree. Although relational structures \cite{peckhaus199919th,tarski1941calculus}, higher order matrices \cite{cayley1894collected} and n-ary operations \cite{peirce1902logic} have been known since the inception of modern mathematics more than a century ago, non-binary instances of such structures have received much less attention than their binary counterparts. More recently, the development of higher category theory \cite{lurie2009higher,simpson2011homotopy}, operads \cite{markl2002operads,leinster2004higher}, opetopes \cite{baez1998higher,cheng2004weak}, hypercompositional algebra \cite{davvaz2009generalization,massouros2021overview} and hyperstructures \cite{baas2019mathematics,baas2019philosophy}, shows a growing interest in ideas that run parallel to the notion of higher arity, despite none of these approaches fully embracing it. In order to capture the essential features of higher arity we shall pay particular attention to ternarity -- the next step above binarity in the arity ladder -- as even the simplest cases display ample perplexing behaviour. Associativity (or lack of thereof) in ternary compositions and ternary products of 3-matrices, adjacency in 3-uniform hypergraphs, and the integration problem for 3-Lie algebras, are all instances where existing binary mathematical formalisms fail to adequately capture the relevant properties of ternary objects. Our observations below strongly suggest the need to develop new mathematical formalisms to faithfully encode genuinely ternary (and higher) structures.\newline

All this compels us to claim that: \textbf{technical and conceptual challenges aside, the potential for scientific discovery and the unparalleled opportunity for mathematical creativity make higher arity research an almost irresistible proposition}.\newline

The outline of this article is as follows: we start with a conceptual definition of `arity' in \ref{arity}; we discuss the ubiquity of binary structures and the need to search for genuine higher arity in \ref{binary}; we analyze a series of ternary mathematical structures and argue for the development of new ternary mathematical theories in \ref{ternarymaths}; we present multiple phenomena from several scientific, social and artistic disciplines to illustrate the ubiquity of higher order phenomena in \ref{ternaryscience}; lastly, we comment on potential future research directions and emphasize the transdisciplinary nature of higher arity science in \ref{concl}.

\section{What is Arity?} \label{arity}

The term `arity' \cite{wiki2022arity} is a noun derived from words such as `binary', `ternary', `n-ary', etc., typically used to describe the number of elements involved in a relation or the number of arguments of an operation. Other terms such as `order', `degree', `adicity', `valency', `rank' or `dimension' are sometimes used to refer to the same concept. However, given the myriad other meanings and connotations that these terms carry across the mathematical sciences, here we propose to fix the meaning of the far less frequently used term `arity'. Still on the topic of nomenclature, we shall refer to the fundamental notion of perceived countable amount as `numerosity' \cite{nunez2017there} to distinguish it from the more technically loaded `quantity' or `cardinality'.\newline

We propose a concise conceptual definition of \textbf{arity} as follows:
\begin{definition}
\textbf{Arity} is the exact numerosity of a relational or functional interdependence between a countable collection of separate entities or discernible parts of a system.
\end{definition}

\begin{figure}[h]
\centering
\includegraphics[scale=0.5]{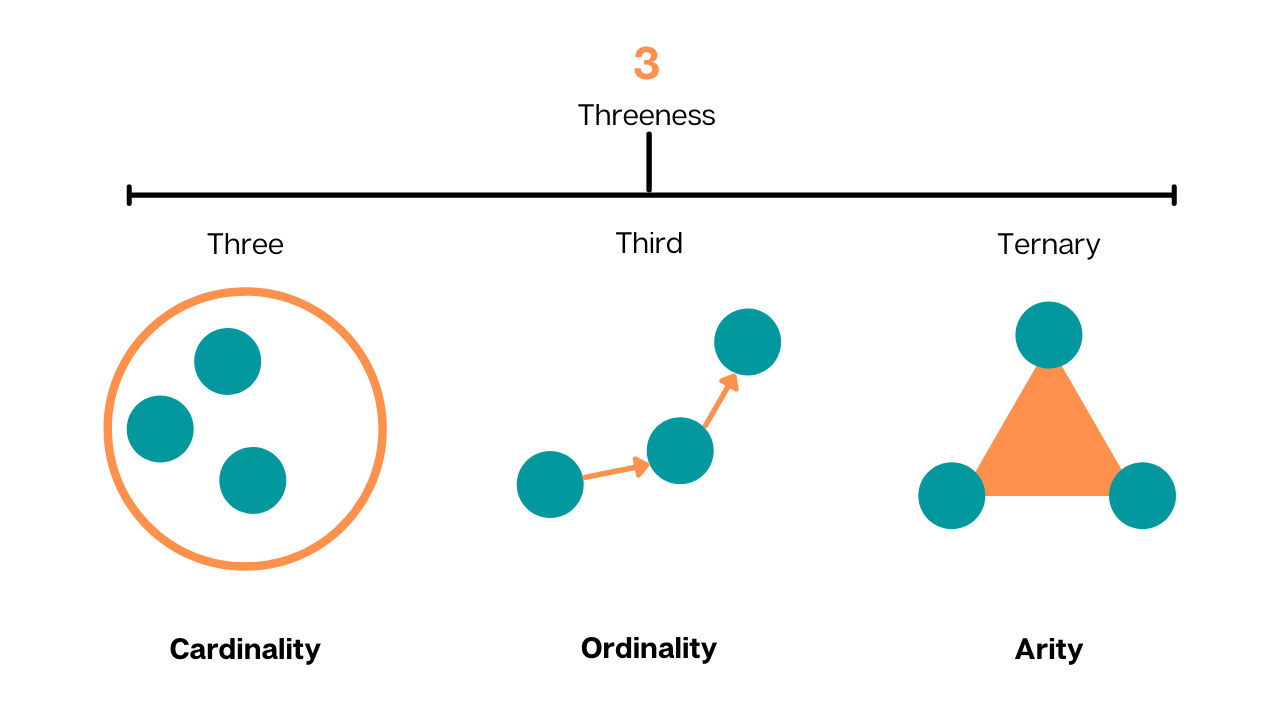}
\caption{The numerosity 3 (orange) in multiple guises for a collection of three objects (teal). Cardinality, perhaps the most direct expression of numerosity, refers to notions of quantity  (comparisons of a set with its parts) or size (comparisons of a set with other sets). Ordinality refers notions of sequence, time-evolution or process between the objects. Arity refers to relations or operations involving all three objects simultaneously.}
\label{VS}
\end{figure}

Arity may be generally understood as an elementary property of an interaction between parts of a system quantifying the amount of discrete components involved in the interaction. In this sense, arity appears as a rudimentary measure of atomic complexity. An effective way to isolate the notion of arity from neighbouring basic concepts related to numerosity is to counterpoise it with cardinality (quantity, size) and ordinality (order, sequence). Figure \ref{VS} illustrates the difference between the concepts of `three', `third' and `ternary'.

\section{The Binarity Bias and Irreducible Higher Arity} \label{binary}

Beyond the trivial arity 1, i.e. the concept of oneness, current intellectual discourse is dominated by one kind of arity alone: binarity. As exemplified by the following cases, the vast majority of formal devices frequently used to articulate models and theories have a very strong binary flavour:
\begin{itemize}
    \item   ``$0+1$'' Operations, ranging from elementary arithmetic to sophisticated notions in abstract algebra, predominantly take two arguments.
    \item   ``$x \sim y$'' Relations are commonly taken between pairs of objects; crucially, equivalence relations and symbolic equality ``='' are defined to be binary.
    \item   ``$A\to B$'' Processes, connections and transformations, abstracted in functions, graphs and categories, rely on input-output or source-and-target paradigms for pairs of objects.
    \item   ``\emph{cogito ergo sum}'' Natural language and truth-valued logic are primarily built from unary assignments or binary operations between sentences and propositions, e.g. subject-verb-object constructions and syllogisms. 
\end{itemize}
Without a broader understanding of the neurobiological and cultural basis of arity, one can only speculate as to what causes this ubiquity of binary structures. We can, nonetheless, identify some factors that are likely to contribute to the status quo:
\begin{itemize}
    \item \textbf{Biological and Cognitive Predispositions:} Some notable features of human biology appear compatible with a preference for binary ideas: bilateral body plan, the perception of chirality (left-right distinction) or sexual reproduction in a binary mating pattern are some examples. Experience of time, causal thinking, perceived sameness of stimuli and limited working memory are aspects of human cognition that also seem to favour binary structures.
    \item \textbf{Symbolic Language Constraints:} Most recorded human languages have a strong sequential nature and tend to string symbols in linear arrays. This is particularly true for the Indo-European languages that found their way into the modern standards of science and mathematics. The development of pre-digital printing technology only exacerbated the preponderance of linear written languages among early scholars as it rendered the efficient communication of non-sequential symbolic information practically impossible. Only recently, with the advent of digital environments, diagrammatic languages and other forms of non-sequential symbolic communication started gaining some traction.
    \item \textbf{Historical Bias:} Beyond just language, the influence of idiosyncratic cultural evolution stretching millennia into the past can be easily recognized in present-day conventions and customs. Modern science and mathematics are relatively recent inventions and, despite the large volume of research output that already exists, only a minuscule fraction of the virtual universe of all possible ideas has been charted. Given their relative simplicity alone, binary structures are expected to be studied extensively before other higher arity generalizations are explored.
\end{itemize}

The thesis that most of current science and mathematics is dominated by binary structures could be contested on the basis of the existence of many examples of operations and relations effectively involving large amounts of elements, e.g. large data sets, multivariable functions, complex algorithms, etc. Although coarse descriptions of such examples will indeed lead to instances of higher arity, a closer look would likely reveal a composite arrangement of fundamentally binary structures. This leads us to the contrasting notions of \textbf{irreducible arity} and \textbf{composite arity}.\newline

To articulate these notions we introduce the concept of \textbf{ariton}\footnote{This term was first suggested by I. Altman during a private conversation.}:
\begin{definition}
An \textbf{ariton} is a minimal expression of a fixed arity. More precisely, a collection of $n$ objects together with an elementary relation or operation between them will be called a $n$-ariton; a concrete pair of related objects is a binariton, a concrete trio of related objects is a ternariton, and so on. 
\end{definition}
Let us clarify this definition with a few examples. Take, for instance, a field of numbers $(\mathbb{F},+,\cdot)$ and consider some elements in it $x,y,z\in\mathbb{F}$. The expressions $x+y$ or $y\cdot z$ are binaritons constructed from the \textbf{basic operations} in $\mathbb{F}$. Expressions involving three elements, such as $x+y+z$ or $x\cdot(y+z)$, formally define ternary operations, however we would not call these ternaritons, but rather composite binaritons, since they are constructed by concatenating binary operations. An example from euclidean geometry is the operation that takes three concurring segments $a,b,c$ as arguments and gives the \textbf{volume of the parallelepiped} spanned by them; the object $V (a,b,c)$ is a ternariton since the volume of a region of space is a primitive notion and all three segments are needed in simultaneous conjunction to define the parallelepiped. In contrast, the analogous construction in linear algebra, i.e. the determinant of the three vectors representing the segments in a $3$-dimensional space $\det(\Vec{a},\Vec{b},\Vec{c})$, is not a ternariton, but rather a composite binariton, since it can be written in terms of the dot and cross products $\det(\Vec{a},\Vec{b},\Vec{c})=\Vec{a}\cdot \Vec{b}\times\Vec{c}$ or, choosing a basis, as a polynomial expression of the vector components, which uses the binary addition and multiplication operations of scalars.\newline

\begin{figure}[h]
\centering
\includegraphics[scale=0.29]{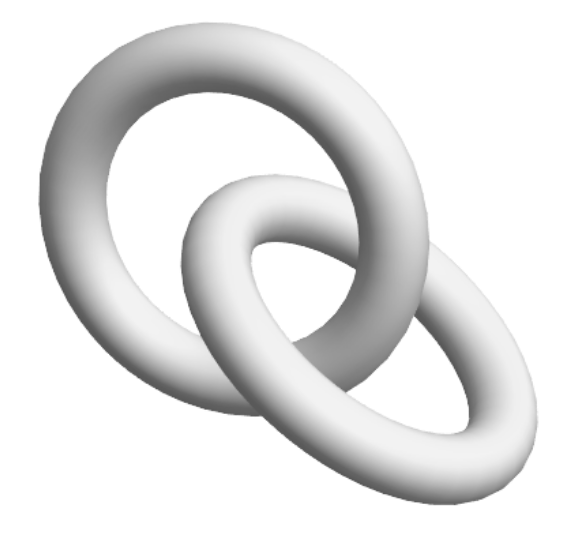}
\caption{Two rings forming a binariton via the linkedness relation.}
\label{linked}
\end{figure}

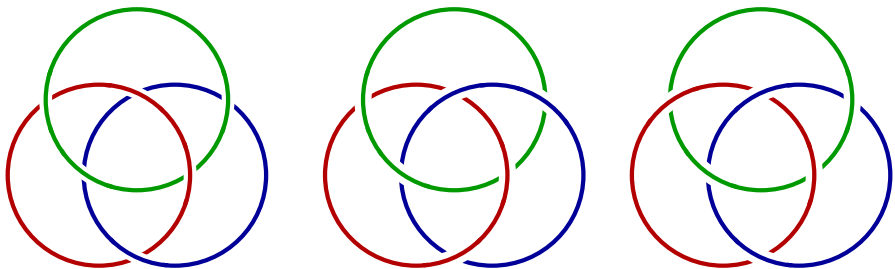
\begin{figure}[h!]
\centering
\begin{tikzpicture}
\begin{knot}[
    clip width=3,
    ]
    \strand [ultra thick, nred  ] (0,0) circle (1.2cm);
    \strand [ultra thick, ngreen] (0.5,1) circle (1.2cm);
    \strand [ultra thick, nblue] (1,0) circle (1.2cm);
    \flipcrossings{2, 4}
\end{knot}
\end{tikzpicture}
\hspace*{0.5cm}
\begin{tikzpicture}
\begin{knot}[
    clip width=4,
    ]
    \strand [ultra thick, nred] (0,0) circle (1.2cm);
    \strand [ultra thick, nblue] (1,0) circle (1.2cm);
    \strand [ultra thick, ngreen] (0.5,1) circle (1.2cm);
    \flipcrossings{1,4,6}
\end{knot}
\end{tikzpicture}\hspace*{0.5cm}
\begin{tikzpicture}
\begin{knot}[
    clip width=4,
    ]
    \strand [ultra thick, nred] (0,0) circle (1.2cm);
    \strand [ultra thick, nblue] (1,0) circle (1.2cm);
    \strand [ultra thick, ngreen] (0.5,1) circle (1.2cm);
    \flipcrossings{1, 2, 5, 6}
\end{knot}
\end{tikzpicture} 
\caption{Knot diagrams for three rings forming (left) an open chain or concatenation of binary links, (centre) a closed chain or triangle of binary links and (right) the Borromean configuration.}
\label{knots}
\end{figure}

\begin{figure}[h!]
\centering
\includegraphics[scale=0.26]{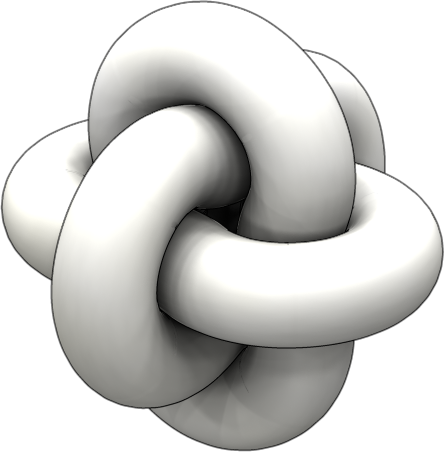}
\caption{Three rings in Borromean configuration forming a ternariton via the linkedness relation.}
\label{borromean}
\end{figure}

The well-known case of the \textbf{Borromean rings} \cite{liang1994borromean,baas2009new} is perhaps the most illustrative example of composite vs irreducible arity. Consider the `linkedness' relation between rings depicted in Figure \ref{linked}. By carefully arranging multiple rings in space, it is easy to see that we can define several linkedness relations of arbitrary arity. Then a natural question arises: given a particular arrangement of rings, are any two of the rings pair-wise linked when considered in isolation from the rest? Knot diagrams help us catalogue topologically different configurations of rings; Figure \ref{knots} shows some possible configurations of three rings based on crossing patterns. The Borromean configuration is remarkable since the isolated pairs of rings are unlinked despite all three being linked in conjunction. In other words, a physical model of the Borromean rings will, in fact, hold together while breaking one of the rings will unlock the other two from each other. This is in stark contrast to the closed chain configuration depicted in Figure \ref{knots}, where breaking one ring does not compromise the linkedness of the remaining two. We thus see that, under the linkedness relation, a chain is a composite binariton while the Borromean configuration is a ternariton (Figure \ref{borromean}).\newline

Aritons appear in a broad spectrum of complex systems as \textbf{atomic cells of interaction}: fundamental forces in multi-particle systems \cite{cohen1993fifty,dobnikar2002many}, processing nodes in computational frameworks \cite{andrews2000foundations,wolfram2021multicomputation}, basic chemical reactions in metabolic networks \cite{wagner2001small,ravasz2002hierarchical}, inter-species relations in ecosystems \cite{schoenly1991trophic,billick1994higher,sole2012self}, protein-protein interactions determining biological function \cite{phizicky1995protein,nooren2003diversity,bertoni2017modeling}, etc. The notion of ariton is implicit in integrated information theory \cite{mediano2022integrated}, aritons occurring across different size scales appear as coarse-graining \cite{flack2017coarse} and aritons persistent in time give rise to the notion of individual \cite{krakauer2020individuality}. From this viewpoint, the general phenomenon of emergence \cite{anderson1972more} may be hypothesized to be a manifestation of arity discrepancies across time and scale that could be measured by some form of arity cohomology.\newline

These first few examples of ternaritons already suggest that the ubiquity of binary structures in the current intellectual mainstream is largely a historical and anthropocentric artifact. This claim will be further consolidated in sections below by the many examples of higher aritons appearing in multiple areas of science and mathematics. 

\section{Ternary Mathematics} \label{ternarymaths}

If we are interested in finding instances of arity higher than 2, it is natural to reach for the next rung in the arity ladder: ternarity. The benefits of focusing on arity 3 are twofold. Firstly, the fact that it is the smallest arity above 2 minimizes the cognitive load when thinking about concrete examples of multi-object systems, as imagining simultaneous interacting components becomes practically unfeasible for human minds beyond numerosities not much larger than 3 \cite{scholl1999tracking}. Secondly, a good understanding of ternarity, in combination with our current knowledge of binarity, will provide a manner of first induction step in the arity ladder similar to how we can think about higher geometric dimensions aided by lower dimensional analogies \cite{blacklock2014analogy}.\newline

The abstract realm of mathematics provides the ideal grounds to systematically investigate ternarity and to get a general sense of its fundamental features. Ternary structures, although rarely found in the existing literature, occupy a natural place in the landscape of mathematical objects and many simple non-trivial examples are readily available. We introduce several such examples in this section by generalizing familiar binary notions such as graphs, matrices, relations, categories and groups. Interestingly, however, making any substantial progress in developing a mathematical theory about these ternary generalizations, one that is at least comparable to those available for their binary counterparts, will turn out to be very challenging. We suspect that among the underlying causes for this difficulty are the limitations imposed by the binarity bias, as explained in \ref{binary}.\newline

A question that can hardly be avoided when thinking about generalizations of binary structures is how to extend the notion of \textbf{sequentiality} into higher arities. A sequence can be understood as iterated binarity, as clearly exemplified by the simple case of a set with a successor function. What are the higher arity analogues of sequences? And, particularly, what is the ternary analogue of a successor function? Although approaches using simplicial complexes \cite{courtney2016generalized} or paths on hypergraphs \cite{carletti2020random} may offer some preliminary answers, there is little consensus on the full extent of these questions.\newline

A particularly profound and far-reaching manifestation of the prevalence of binary structures in mathematics is the ubiquity of \textbf{associativity}. From transitivity of relations in sets, to operations in algebras, to sequences of processes in algorithms, to composition of morphisms in categories, to concatenation of paths in networks, associativity is deeply ingrained in the vast majority of mathematical formalisms in use today. Given how universal associativity seems to be for binarity, one may hypothesize that a similar property plays an analogous role for ternarity. Could a form of generalized associativity, such as the na\"ive ternary extension of binary associativity
\begin{equation*}
    (abc)de = a(bcd)e = ab(cde),
\end{equation*}
have a similarly privileged place at the core of ternary mathematics? Excitingly, the answer to this question is not yet known. Ternary associativity is a very elusive concept and it appears that the obvious generalizations of binary associativity are ill-suited to capture the algebraic properties of even the simplest examples of ternary structures. Despite these difficulties, we shall present compelling evidence that points to the existence of a well-defined, if perhaps non-unique, notion of ternary associativity.\newline

In the sections to follow we aim to highlight fundamental structural properties of ternary objects, therefore, we shall prioritize describing their compositional and algebraic features over combinatorial ones. More concretely, this means that we will remain deliberately agnostic with respect to the precise formalization of higher arity ordering and permutation symmetry, using unordered sets whenever possible and otherwise choosing convenient (possibly non-unique) orderings in the interest of conciseness. A more detailed treatment of these topics will follow in future work by the authors \cite{zapata2022triality}.

\subsection{$3$-Uniform Hypergraphs}\label{graph}

The most minimalist mathematical expression of arity is encapsulated in the notion of hypergraph. Ordinary graphs, which we take as $2$-uniform simple hypergraphs on a set of vertices \cite{voloshin2009introduction,bretto2013hypergraph}, encode information in the form of links between pairs of vertices and we thus regard them as binary structures. A natural ternary generalization is to consider $3$-uniform simple hypergraphs, which encode information in the form of ternary links between trios of vertices. Research on hypergraphs typically focuses on general traits and large scale structure, as seen in recent work on combinatorics \cite{lu2016connected,conlon2019hypergraph,veldt2020hypergraph}, complex networks \cite{battiston2021physics,neuhauser2021consensus}, higher order dynamics \cite{carletti2020random,carletti2020dynamical,landry2021hypergraph} or computational complexity \cite{leordeanu2012efficient,lanzinger2020hypergraph}. In contrast, our interest in hypergraphs is much more humble as we will simply need to discuss some elementary topological properties of small uniform hypergraphs.\newline

A connected $k$-uniform simple hypergraph will be called a $k$\textbf{-graph} for short. We shall regard hypergraphs as a general template to study the basic compositional properties of $n$-ary relations and operations. In this view, the emphasis is placed on hyperedges and their connectivity instead of the underlying vertices -- a similar perspective to how category theory emphasizes morphisms and their compositions over individual objects \cite{leinster2014basic}. Accordingly, we define the \textbf{size} of a $k$-graph simply as the total number of edges, i.e. the cardinality of the set of edges. For any given $k$-graph we distinguish two kinds of vertices according to their degree, i.e. the number of edges containing them: vertices of minimum degree are called \textbf{external vertices} and all other vertices are called \textbf{internal vertices}, if all the vertices have the same degree they are all defined to be internal except when the $k$-graph is of size $1$ (the unique smallest $k$-graph), in which case they are defined to be external. The number of external vertices of a $k$-graph will be called its \textbf{externality}.\newline

\begin{figure}[h]
\centering
\includegraphics[scale=0.37]{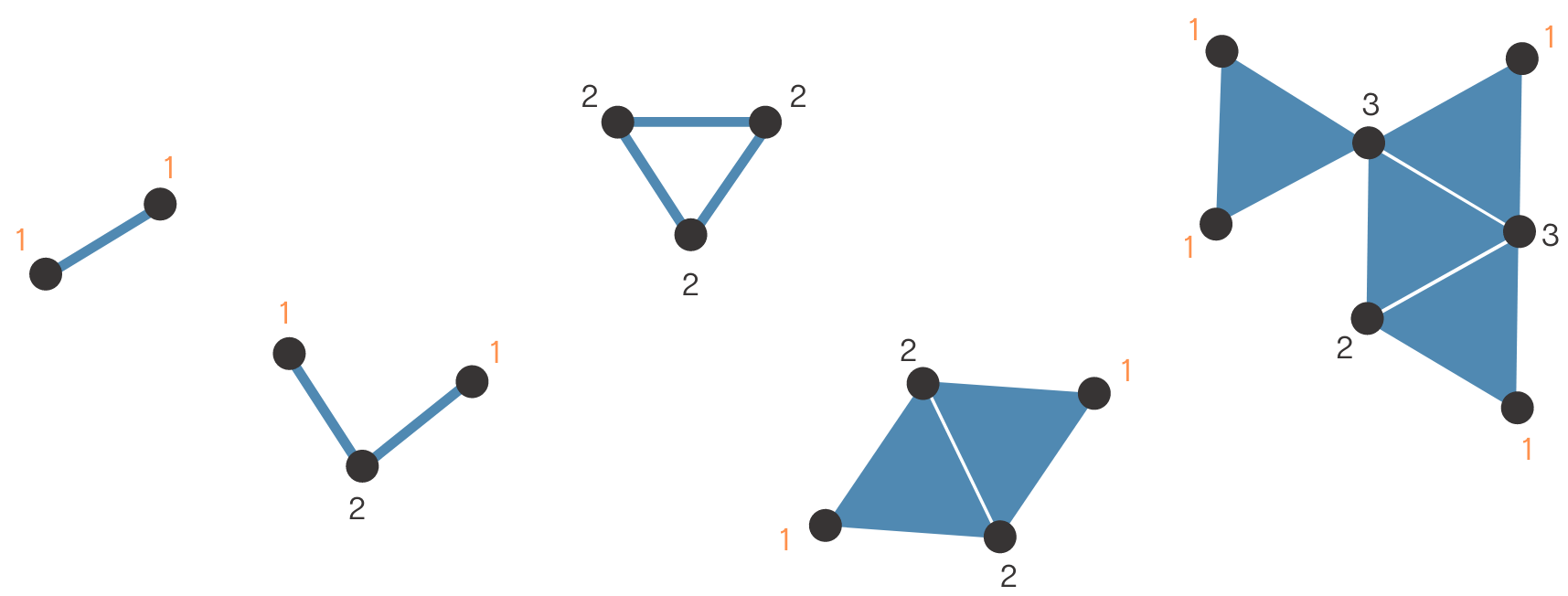}
\caption{Some $2$-graphs and $3$-graphs displaying degree counts, orange highlights correspond to the minimum degrees defining the external vertices.}
\label{graph1}
\end{figure}

The central notion that we would like to put forward in this section is what we call \textbf{motif adjacency}. To motivate this idea let us first reinterpret the familiar notion of adjacency in ordinary graphs. Given a particular graph $G$, adjacency is simply captured as a Boolean function on pairs of vertices defined to be $1$ when there is an edge between them and $0$ otherwise, the information of this Boolean function is encoded in the so-called adjacency matrix $A_G$. It is well known that powers of the adjacency matrix $A^l_G$ encode information about the existence of paths of length $l$ between pairs of vertices. Conversely, keeping the matrix algebra on Boolean values, each higher power of the adjacency matrix encodes the data of a graph that we denote by $G^l$. The graph $G^l$ has the same vertex set as $G$ and an edge between a pair of vertices of $G^l$ is drawn if and only if there exists at least one path of length $l$ in the original graph $G$ connecting them. Our key observation is that the graph $G^l$ can be generated by a motif detection algorithm: we choose our motif to be the graph $P_{l+1}$, consisting of a chain of $l$ edges with an external vertex at each end, and define an edge of $G^l$ between a pair of vertices whenever there exists a subgraph of $G$ that is isomorphic to the motif $P_{l+1}$ and whose external vertices coincide with the pair of vertices under consideration. A pair of vertices of $G$ that are linked by an edge in $G^l$ are said to be \textbf{motif-adjacent} or, specifically, $P_{l+1}$\textbf{-adjacent}.\newline

Although anyone familiar with ordinary graph theory would probably find this characterization of adjacency unnecessarily convoluted, our approach opens up a straightforward way to define adjacency in general hypergraphs: consider a hypergraph $H$ and a $k$-graph $M$ of externality $m\geq 2$ to be taken as the motif, then, a collection of vertices $\{i_1,i_2,\dots,i_m\}$ of $H$ is said to be $M$\textbf{-adjacent} if there exists a subhypergraph $S\subset H$ whose external vertices are precisely $\{i_1,i_2,\dots,i_m\}$ and such that $S\cong M$. We define the $M$\textbf{-adjacency hypergraph derived from} $H$ by the motif detection algorithm described above in the obvious way. Note that the $M$-adjacency hypergraph is $m$-regular by construction. A thorough treatment of motif adjacency, involving multiway systems of hypergraph rewrites based on adjacency relations, would take us too far beyond the scope of the present article. The authors will revisit this topic in connection with ternary associativity in future work \cite{zapata2022heap}.\newline

\begin{figure}[h]
\centering
\includegraphics[scale=0.37]{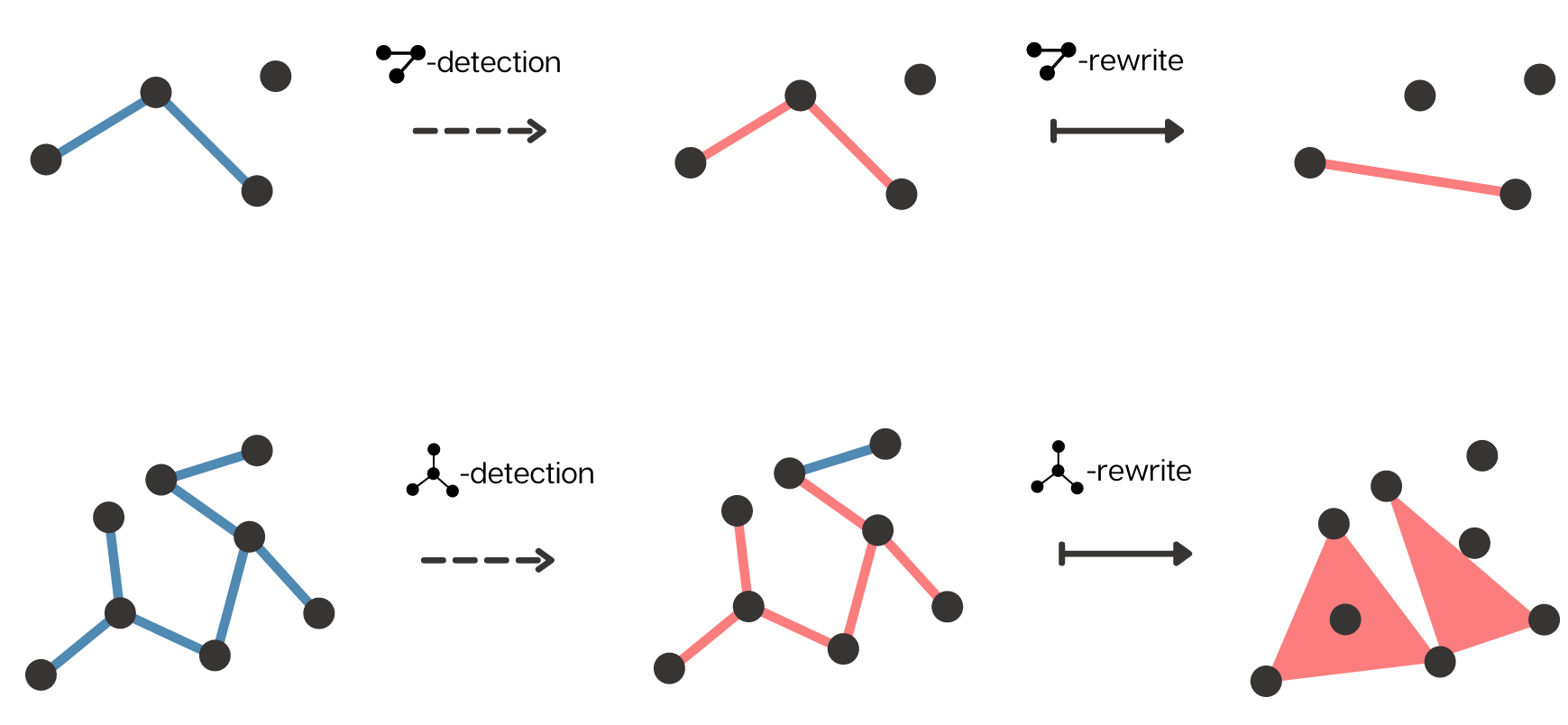}
\caption{Adjacency rewrites.}
\label{graph2}
\end{figure}

Under the lens of motif adjacency, the conventional treatment of adjacency in ordinary graphs appears somewhat accidental in that the motifs chosen ($P_l$) always happen to have two external vertices, thus making adjacency-derived hypergraphs always $2$-regular, that is, ordinary graphs again. This may give the false impression that adjacency rewrites are internal operations in the class of $2$-graphs, however, as illustrated in Figure \ref{graph2}, this is not the case when more general motifs are considered. Our notion of motif adjacency can be thought of as a way to allow for more topological diversity in the probing tools to investigate graph structure: all $P_l$ graphs are topologically equivalent (homeomorphic to a segment regarded as topological spaces) while equally simple motifs, such as the claw $K_{1,3}$ \cite{graphclass2021claw} in Figure \ref{graph2} or small graphs containing cycles like the bull \cite{graphclass2021bull}, display different topologies. This being said, there is good reason for the extended use of $P_l$ as motif in ordinary graphs: a $P_l$ motif corresponds to the notion of path of length $l$, arguably the most central concept in graph theory.\newline

Recall that our goal was to extend familiar ideas of graph theory from $2$-graphs to $3$-graphs. Although it may be debatable whether it is worthwhile to consider general motif adjacency in $2$-graphs, the reality is that when it comes to $3$-graphs we do not have much of a choice. This is directly connected with the sequentiality problem: there is no consensus on what constitutes a satisfactory higher order generalization of a path or a sequence. Even if we accept some of the definitions of paths in hypergraphs that have been proposed \cite{gallo1993directed,luczak2019paths,carletti2020random}, we must reconcile the fact that ternary (and higher) sequentiality is no longer uniquely characterized.\newline

\begin{figure}[h]
\centering
\includegraphics[scale=0.42]{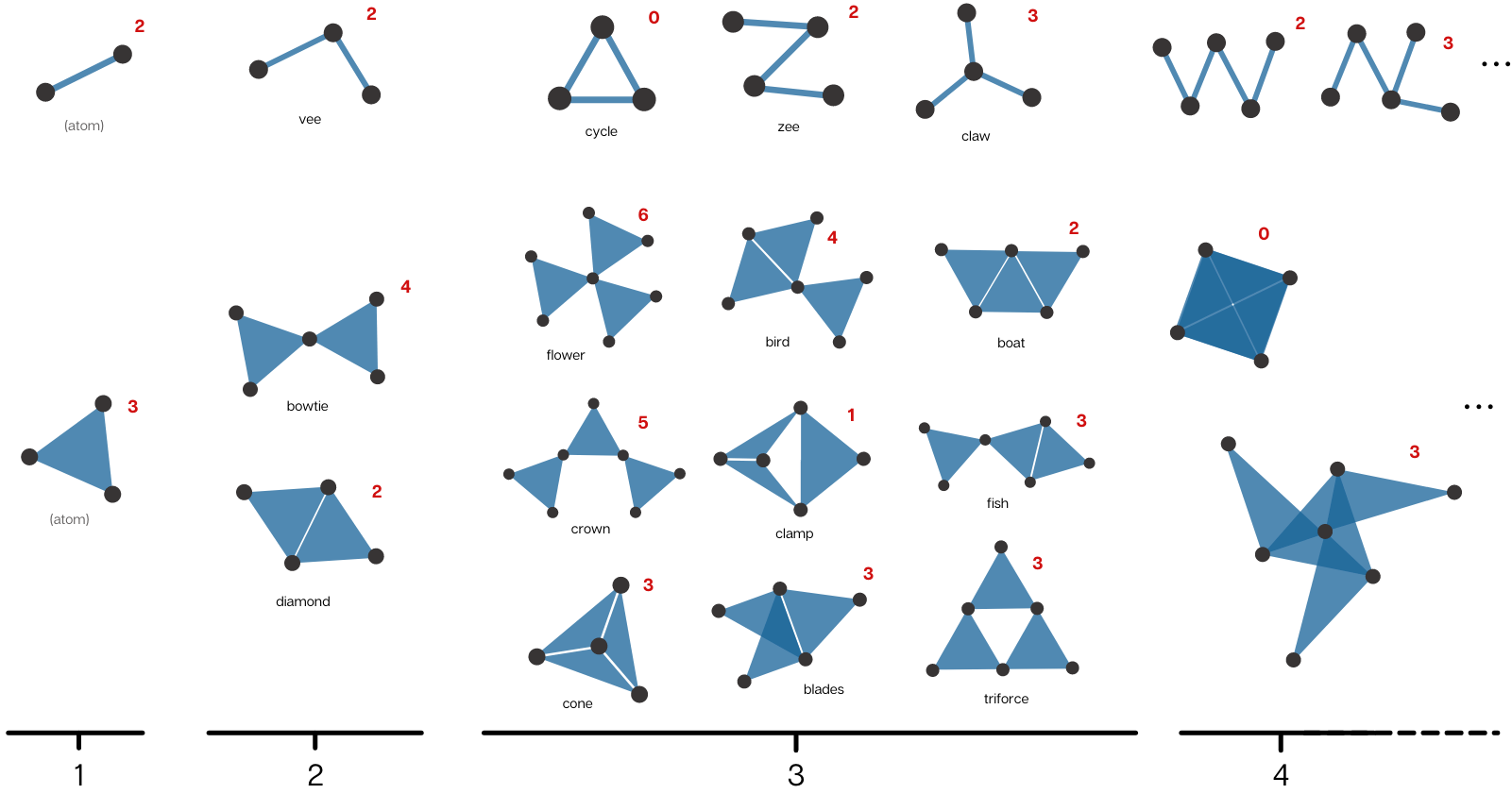}
\caption{Isomorphism classes of $2$-graphs and $3$-graphs of increasing size labelled by externality.}
\label{graph3}
\end{figure}

A parsimonious and unbiased approach would suggest that we consider all possible isomorphism classes of small $k$-graphs as potential motifs for adjacency rewrites. To this end, we catalogue small $2$-graphs and $3$-graphs according to their size and externality in Figure \ref{graph3} giving them short suggestive names for later reference. We use a hypergraph generalization of the NAUTY graph isomorphism algorithm \cite{mckay1981practical,mckay2014practical} implemented in the SageMath software \cite{sage2022nauty} to generate all the isomorphism classes of hypergraphs of a fixed size. If we simply count the number of isomorphism classes of $2$-graphs and $3$-graphs of increasing size we find a first instance of the stark contrast between binarity and ternarity:

\begin{center}
\begin{tabular}[c]{ |c|c|c|c|c|c|c|c|c|c| } 
 \hline
 size & 1 & 2 & 3 & 4 & 5 & 6 & 7 & 8 & 9 \\
 \hline
 \# $2$-graph isomorphism classes & 1 & 1 & 3 & 5 & 12 & 30 & 79 & 227 & 710 \\
 \hline
 \# $3$-graph isomorphism classes & 1 & 2 & 9 & 51 & 361 & 3683 & 47853 & 780725 & 15338226 \\ 
 \hline
\end{tabular}
\end{center}

Similar to how adjacency in ordinary graphs can be understood as iterated vee motif detection, i.e. paths are sequences of edges, it is natural to first consider the smallest $3$-graphs of Figure \ref{graph3} as motifs and then build larger motifs by iteration. In this sense, given a collection of $k$-graphs, we define their \textbf{splicings} as new $k$-graphs resulting from attachments of external vertices. It is easy to see that $2$-graphs and $3$-graphs in Figure \ref{graph3} are generated in this way from splicings of smaller graphs. Also note how splicings recombine externalities in multiple ways: from the atomic $2$-graph of externality 2 and the atomic $3$-graph of externality 3 we find splicings of externalities 0, 1, 4, 5... This motivates the notion of \textbf{externality-preserving splicings} or \textbf{compositions}, defined as any splicing of a collection of $k$-graphs of the same externality $n$ that results in a $k$-graph of externality $n$. The smallest composition in the class of $2$-graphs is given by the vee shape and its iteration recovers the usual notion of paths in ordinary graph theory. In contrast, the smallest composition in the class of $3$-graphs is no longer unique: the cone, blades, fish and triforce shapes all represent splicings preserving externality 3.\newline

The following sections \ref{cubix}, \ref{rel} and \ref{cat} will vindicate the compositional framework of $k$-graphs described in this section as an effective notational paradigm to capture the basic features of higher order matrix multiplications, higher arity relation composition and higher order isomorphism compositions, respectively. Particularly, we will see how the cone, the blades, the triforce and the fish shapes appear naturally as canonical ternary compositional patterns in all those different contexts.

\subsection{$3$-Matrix Algebras}\label{cubix}

Consider some field\footnote{We could take a ring more generally since we really only make use of the additive and multiplicative structures, existence of multiplicative inverses or commutativity play no role in our discussion.} of numbers $(\mathbb{F},+,\cdot)$. We define a $n$\textbf{-matrix} as a multi-index array of the form $[a_{i_1i_2\cdots i_n}]$, where $a_{i_1i_2\cdots i_n}\in \mathbb{F}$ and each of the $n$ indices ranges in some finite interval $i_j\in \{1,2,\dots,N_j\}$, the total number of indices $n$ is called the \textbf{dimension} of $a$ and the tuple of lengths of the index ranges $(N_1,N_2,\dots,N_n)$ is called the \textbf{size} of $a$. Although known to the pioneers of matrix theory J. Sylvester and A. Cayley as far back as the mid 19$^{\text{th}}$ century \cite{cayley1894collected}, such higher order matrices have been largely ignored during the history of modern mathematics. In recent years, however, $n$-matrices have received renewed interest from multiple communities: they appear in modern revisions of the classical theory of discriminants \cite{gelfand2008discriminants}; they offer natural examples of $n$-ary algebras \cite{abramov2009algebras,azcarraga2010nary,abramov2011ternary,bai20143lie}, particularly, in applications in mathematical physics \cite{kerner1997cubic,kerner2008ternary}; they encode the adjacency information of regular hypergraphs \cite{gnang2014combinatorial,courtney2016generalized,battiston2020networks}; they capture the algebraic data of association schemes \cite{mesner1990association,mesner1994ternary}; computer scientists use them as data structures \cite{shen2019simplified,ji2019survey}; and they provide the natural grounds to explore higher order generalizations of the eigenvalue problem \cite{gnang2011spectral}, invertibility \cite{gnang2020bhattacharya} or canonical forms \cite{gnang2020symmetrization}.\newline

Our goal is to investigate the basic algebraic structure of $n$-matrices but before we begin our discussion we must warn against the temptation to consider $n$-matrices as tensors in the ordinary sense of linear algebra. This is, in fact, a common practice -- and the source of a great deal of confusion -- in some of the research areas mentioned above, particularly in machine learning \cite{bratseth2021tensor} and physics \cite{kerner1997cubic}. It is important to keep in mind that the basis expression of a tensor in conventional linear algebra is not just a multi-index array of scalars but an entire equivalence class of such arrays defined by the covariance relation induced by all the possible choices of basis. We emphasize that we take $n$-matrices simply as arrays of scalars, without any further structure, and, although they will be shown to share many features with ordinary (rectangular) matrices, we shall remain open-minded and let them speak for themselves before we bin them in an existing mathematical formalism.\newline

A $n$-matrix $a=[a_{i_1i_2\cdots i_n}]$ can be effectively thought of as an element of the set $\mathbb{F}^{N_1\times N_2 \times \cdots \times N_n}$. This suggests a hypergraphic notational convention where we place the $n$-matrix $a$ as a hyperedge on vertices $\{\mathbb{F}^{N_1}, \mathbb{F}^{N_2}, \dots, \mathbb{F}^{N_n}\}$. As can be seen in the many applications mentioned above, the most commonly occurring elementary $n$-matrix operations are what we call \textbf{index splicings}. Given a collection of matrices of generally different dimensions and sizes, index splicings are constructions that produce new matrices provided some of the indices range over the same interval. There are two kinds of splicing operations: by \textbf{incidence}, defined simply as multiplication of matrix elements along the common indices, and by \textbf{contraction}, defined by \emph{summing over} or \emph{tracing out} the common indices. Let us illustrate these. Take a $n$-matrix $a_{i_1i_2\cdots i_n}$, a $m$-matrix $b_{j_1j_2\cdots j_m}$, a $l$-matrix $c_{k_1k_2\cdots k_l}$ and suppose that the indices $i_1$, $j_2$ and $k_l$ all range over the same interval $\{1,2,\dots,N\}$, an incidence splicing of $a$ and $b$ results in a $(n+m-1)$-matrix $r$ whose entries are given by
\begin{equation*}
    r_{\textcolor{nred}{k}i_2\cdots i_nj_1\cdots j_m}= a_{\textcolor{nred}{k}i_2\cdots i_n} \cdot b_{j_1\textcolor{nred}{k}\cdots j_m},
\end{equation*}
a contraction splicing of $a$, $b$ and $c$ results in a $(m+n+l-3)$-matrix $s$ whose entries are given by
\begin{equation*}
    s_{i_2\cdots i_nj_1\cdots j_mk_1k_2\cdots k_{l-1}}=\sum_{\textcolor{nred}{p}=1}^N a_{\textcolor{nred}{p}i_2\cdots i_n} \cdot b_{j_1\textcolor{nred}{p}\cdots j_m} \cdot c_{k_1k_2\cdots \textcolor{nred}{p}}.
\end{equation*}
Note that general splicing constructions may combine several simultaneous incidences and contractions of matrices. Index splicings are natural minimal operations that only exploit the array structure of matrices and the additive and multiplicative operations of their scalar entries. Since splicings can only be defined on indices of the same size, the hypergraph notation suggested above naturally captures index splicing operations as $n$-graph splicings in the sense defined in \ref{graph}. Let us illustrate the \textbf{matrix splicing notation} by applying it to some elementary constructions with $2$-matrices. Given two $2$-matrices $a$, $b$ sharing an index of common size, we can define their ordinary associative multiplication by a contraction splicing or we can define a $3$-matrix by an incidence splicing:
\begin{center}
    \includegraphics[scale=0.2]{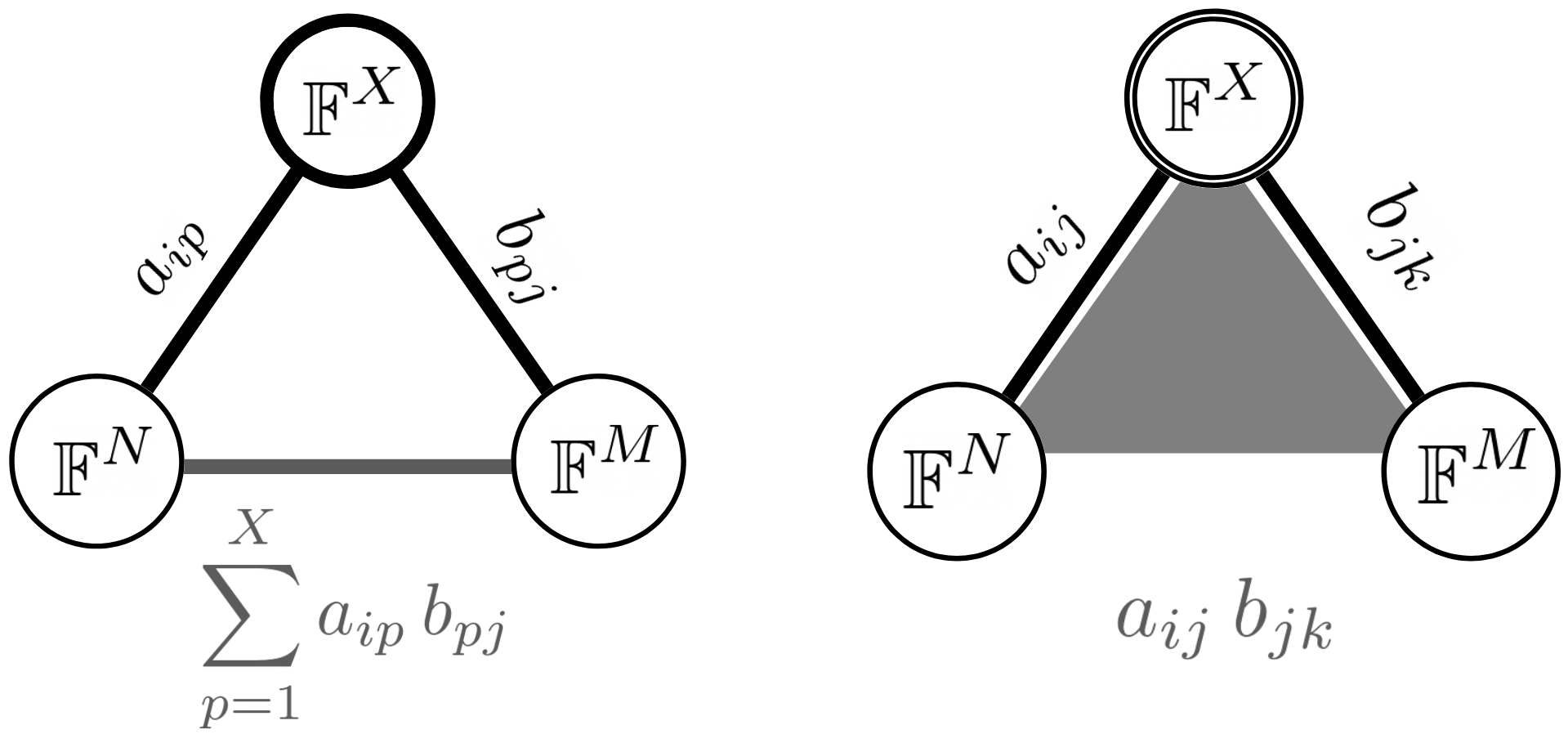}
\end{center}
where both splicings correspond to the same underlying $2$-graph shape, the vee. For three $2$-matrices $a$, $b$, $c$ more $2$-graph shapes are available to define splicings. Here are two splicings defined with the cycle and claw shapes:
\begin{center}
    \includegraphics[scale=0.24]{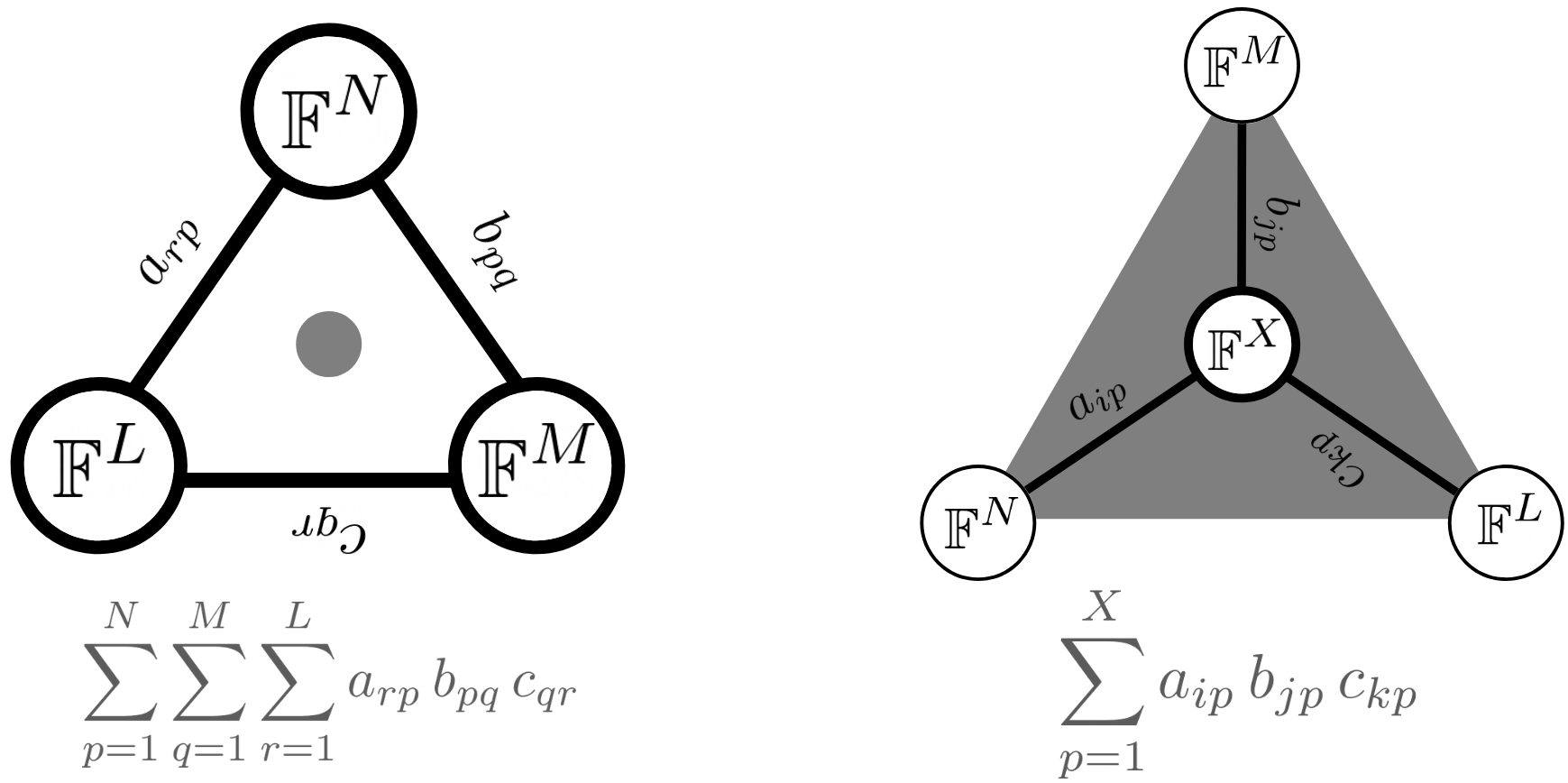}
\end{center}

We recognize the familiar multiplication and trace operations for ordinary (rectangular) matrices in the splicings above but we also find two operations that combine $2$-matrices to give a $3$-matrix; this illustrates that splicings of matrices generally mix dimensions. The hypergraph notation we have introduced allows us to see that an effort to catalogue matrix splicings of increasing dimension would be analogous to our classification of the isomorphism classes of small $n$-graphs in \ref{graph}. If we account for two kinds of vertices corresponding to incidence and contraction splicings, finding \textbf{dimension-preserving splicings}, like the familiar case of $2$-matrix multiplication, becomes a simple combinatorial task on the vertices of small isomorphism classes of hypergraphs.\newline

According to Figure \ref{graph3}, the next higher instances of dimension-preserving splicings occur for a pair of $3$-matrices in diamond configuration and for trios of $3$-matrices in cone, blades, triforce, clamp, boat and fish configurations. Let us define some of these explicitly by making use of our matrix splicing notation. Given two $3$-matrices $a$, $b$ sharing two indices of common size, we can arrange them in a diamond shape with an incidence and a contraction to yield a $3$-matrix again:
\begin{center}
    \includegraphics[scale=0.2]{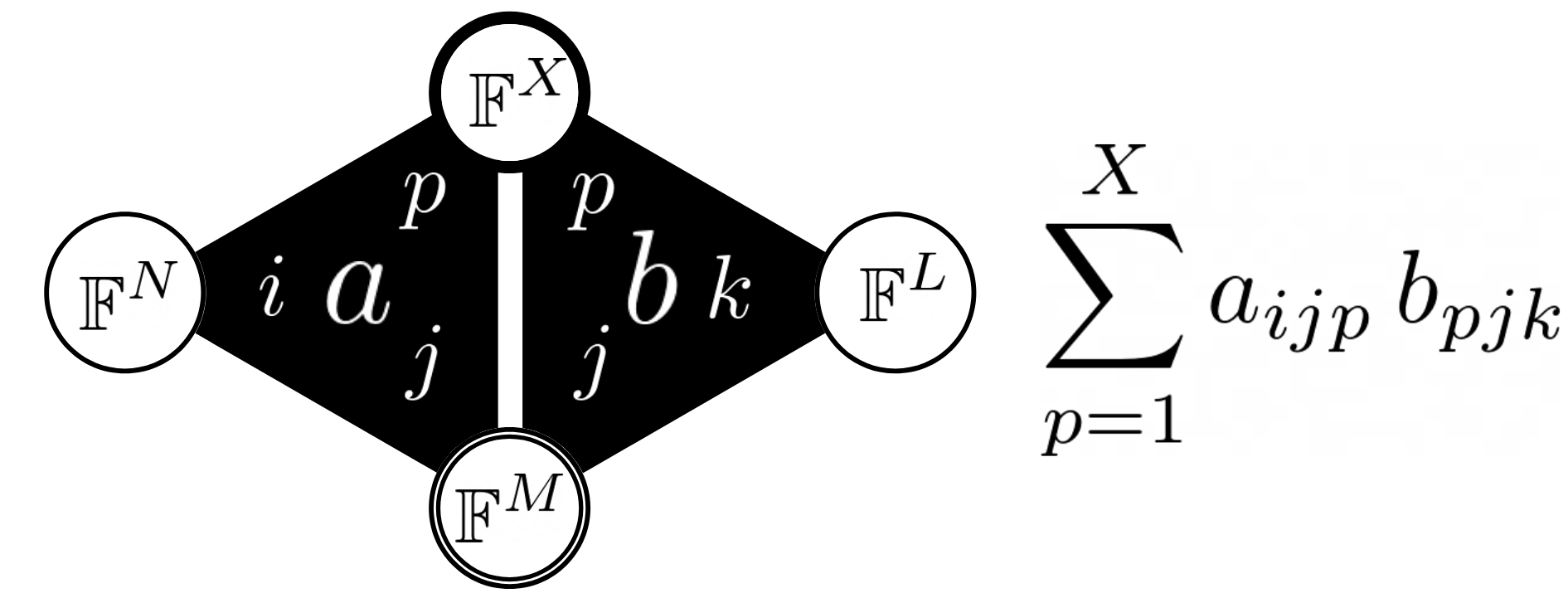}
\end{center}
Now take three $3$-matrices $a$, $b$, $c$. Using the cone shape with a contraction in the index shared by the three $3$-matrices and incidences in the other pairs of shared indices we define the following operation:
\begin{center}
    \includegraphics[scale=0.2]{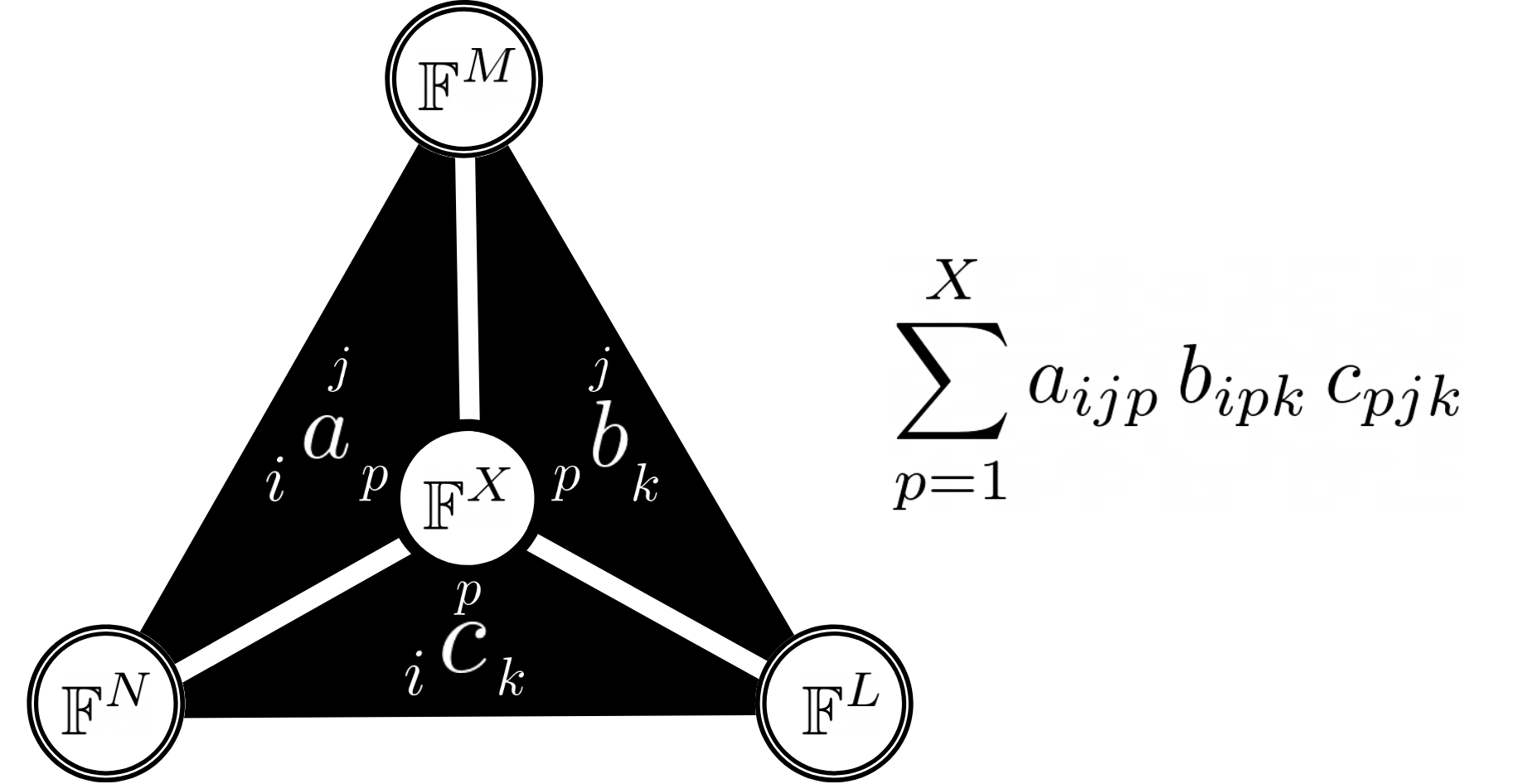}
\end{center}
Using the triforce shape with contractions over all shared indices we find:
\begin{center}
    \includegraphics[scale=0.2]{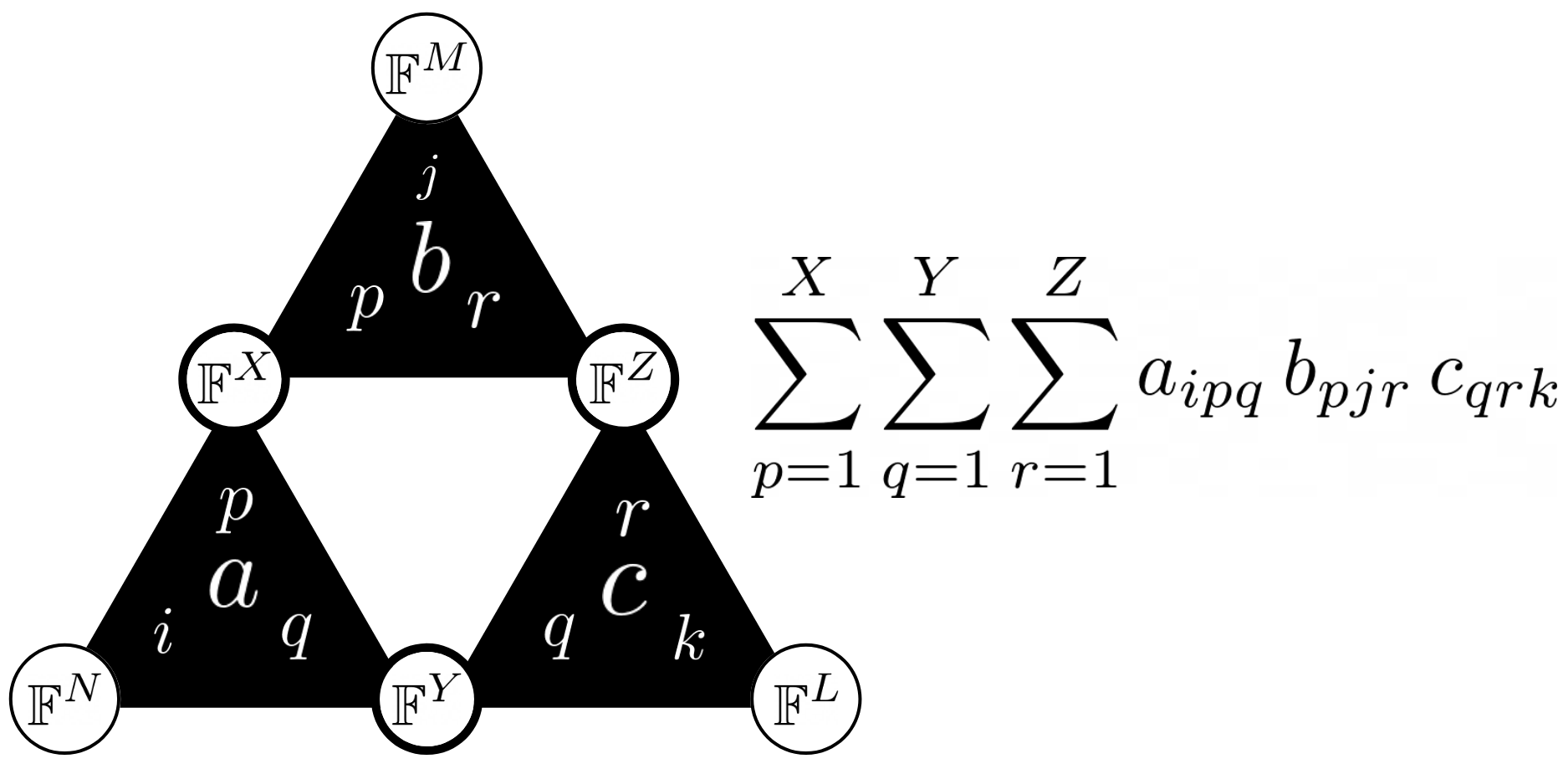}
\end{center}
As we anticipated in our note about ternary order at the start of this section \ref{ternarymaths}, orderings are left deliberately loose to emphasize compositional patterns over combinatorial ones. For instance, by fixing a sequential order of indices, the formulas we have given become well-defined in terms of standard multi-index array notation, however, in doing so, we are making non-canonical choices and multiple order-inequivalent splicings can then be defined by combinatorial permutations of indices. We shall not discuss such order-inequivalent splicings here.\newline

Dimension-preserving splicings induce internal operations for classes of \textbf{regular} $n$-matrices, i.e. $n$-matrices $[a_{i_1i_2\cdots i_n}]$ whose indices all range over the same interval $i_j\in \{1,2,\dots,N\}$. We call a dimension-preserving splicing operation a \textbf{multiplication}. The lowest-dimensional instance of such an internal operation is the single-index contraction in the class of regular $2$-matrices, that is, the familiar unital associative algebra of square matrices. Going one arity step higher, we see from the examples above that splicings give ternary multiplications for regular $3$-matrices. We shall devote the rest of this section to the study of the basic properties of these ternary algebras.\newline

A regular $3$-matrix is often called a `cubic matrix' in the existing literature \cite{kerner1997cubic,kerner2008ternary,abramov2011ternary}. In the interest of brevity, we propose the portmanteau `\textbf{cubix}' to refer to the same concept. The \textbf{space of cubices} of size $N$ over a field (or ring) $\mathbb{F}$, denoted by $\mathbb{F}^{N^3}$, is defined as the set of regular $3$-matrices $a=[a_{ijk}]$, where $i,j,k\in \{1,2,\dots,N\}$, endowed with the obvious entry-wise $\mathbb{F}$-vector space structure. By considering the $3$-graphs of size 3 from Figure \ref{graph3} and all the topologically-distinct distributions of index incidences and contractions, we can define the following \textbf{inequivalent ternary multiplications} for any given trio of cubices $a,b,c\in\mathbb{F}^{N^3}$:
\begin{align*}
    p_1(a,b,c)_{ijk}&:= \sum_{p=1}^N a_{ijp}\, b_{ipk}\, c_{pjk}  &\text{(cone)}\\
    p_2(a,b,c)_{ijk}&:= \sum_{p,q=1}^N a_{ipq}\, b_{pjq}\, c_{pqk} &\text{(blades)}\\
    p_3(a,b,c)_{ijk}&:= \sum_{p,q,r=1}^N a_{ipq}\, b_{pjr}\, c_{qrk} &\text{(triforce)}
\end{align*}

\begin{align*}
    p_4(a,b,c)_{ijk}&:= \sum_{p,q,r=1}^N a_{ijp}\, b_{pqr}\, c_{qrk} &\text{(fish)}\\
    p_{5.1}(a,b,c)_{ijk}&:= \sum_{p,q=1}^N a_{ijp}\, b_{pjq}\, c_{qjk} &\text{(boat)}\\
    p_{5.2}(a,b,c)_{ijk}&:= \sum_{p,q=1}^N a_{ipq}\, b_{pjq}\, c_{pjk} &\text{(boat)}\\
    p_{6.1}(a,b,c)_{ijk}&:= \sum_{p,q=1}^N a_{ijp}\, b_{ijq}\, c_{pqk} &\text{(clamp)}\\
    p_{6.2}(a,b,c)_{ijk}&:= \sum_{p,q=1}^N a_{ijp}\, b_{ipq}\, c_{qjk} &\text{(clamp)}\\
    p_{6.3}(a,b,c)_{ijk}&:= \sum_{p,q=1}^N a_{ijk}\, b_{pjq}\, c_{pqk} &\text{(clamp)}\\
    p_7(a,b,c)_{ijk}&:= \sum_{p=1}^N a_{ijk}\, b_{ijp}\, c_{pjk}  &\text{(cone)}
\end{align*}

Some of these have been identified in the literature: $p_1$ was originally found in the context of relational schemes \cite{mesner1990association,mesner1994ternary} and has been called the Bhattacharya-Mesner product in later work that further investigated its properties \cite{gnang2020bhattacharya}, $p_3$ and $p_4$ appeared in early work on ternary relational structures \cite{mcculloch1967triadas} and are both well-known in recent cubic matrix research \cite{abramov2011ternary}. To the best of our knowledge, $p_2$, $p_5$, $p_6$ and $p_7$ are novel and introduced here for the first time.\newline

Each of the ternary multiplications induces a different ternary algebra structure on the space of cubices $(\mathbb{F}^{N^3},p_i)$. Although they all coexist simultaneously on $\mathbb{F}^{N^3}$, we shall consider them separately and postpone for future work the study of their interaction. This is particularly useful to highlight the analogies between ternary algebras of cubices and the binary algebra of square matrices. Note that multiplications $p_1$, $p_2$ and $p_3$ all involve index incidences and contractions of three cubices in symmetrical roles, similar to how multiplication of ordinary matrices involves index contractions of two matrices in symmetrical roles. Therefore, we argue that $p_1$, $p_2$ and $p_3$ are the direct ternary generalizations of the ordinary multiplication of matrices. In contrast, all the other multiplications are asymmetric in their splicing patterns. For instance, $p_4$ involves cubices arranged in a manner of input-throughput-output scheme or, as pictorially suggested by the fish shape, a tail-body-head scheme. We shall focus on the algebras $(\mathbb{F}^{N^3},p_i)$ $i=1,2,3,4$, which we call the \textbf{cone}, \textbf{blades}, \textbf{triforce} and \textbf{fish cubix algebras}, respectively.\newline

The notion of identity matrix $[\delta_{ij}]$ as the neutral element for ordinary matrix multiplication finds natural analogues in the space of cubices $\mathbb{F}^{N^3}$. We can define the \textbf{partial identities} $i_1$, $i_2$ and $i_3$ as stacks of ordinary identity matrices in each of the three directions of the cubic arrays:
\begin{equation*}
    [i_1]_{ijk}= \delta_{ij}\qquad [i_2]_{ijk}= \delta_{jk}\qquad [i_3]_{ijk}= \delta_{ik}.
\end{equation*}
We can also define the \textbf{tridentity} as the unit diagonal cubix:
\begin{equation*}
    [I]_{ijk}=\delta_{ijk}.
\end{equation*}
These four special cubices $i_1,i_2,i_3,I\in\mathbb{F}^{N^3}$ satisfy identity-like properties with respect to ternary cubix multiplications. With the aid of a simple computer algebra implementation of cubix multiplication in Mathematica software we were able systematically list all such properties. Let an arbitrary $a\in\mathbb{F}^{N^3}$, then we have the identity properties of the cone algebra:
\begin{equation*}
    p_1(a,i_2,i_3)=a \qquad p_1(i_2,a,i_1)=a \qquad p_1(i_3,i_1,a)=a
\end{equation*}
the identity properties of the blades algebra:
\begin{equation*}
    p_2(a,i_1,i_3)=a \qquad p_2(i_1,a,i_3)=a \qquad p_2(i_3,i_1,a)=a
\end{equation*}
the identity properties of the triforce algebra:
\begin{align*}
    p_3(a,i_1,I)=a \qquad p_3(i_1,a,I)=a \qquad p_3(i_3,I,a)=a\\
    p_3(a,I,i_3)=a \qquad p_3(I,a,i_2)=a \qquad p_3(I,i_2,a)=a
\end{align*}
and the identity properties of the fish algebra:
\begin{equation*}
    p_4(a,i_1,I)=a \qquad p_4(a,i_3,I)=a \qquad p_4(a,I,i_2)=a \qquad p_4(a,I,i_3)=a \qquad p_4(a,I,I)=a
\end{equation*}
Furthermore, the four identity cubices $\{i_1,i_2,i_3,I\}$ form partial subalgebras in each of the cubix algebras, in the sense that properties such as $p_1(i_2,i_1,i_3)=I$ or $p_3(I,I,I)=I$ hold. With some possible redundancies due to the general identity equations above, there are 63 such properties in the cone algebra, 52 in the blades algebra, 31 in the triforce algebra and 28 in the fish algebra. All these represent another staggering contrast between binarity and ternarity: what was a single element satisfying a single identity property in the binary algebra of square matrices turns into several identity elements satisfying a plethora of identity-like properties in the ternary algebras of cubices. To further underline the kaleidoscopic nature of the phenomenon of identity in cubix algebras, note that currying identity matrices into cubix multiplications gives binary operations of cubices, for instance $p_1(\cdot,\cdot,i_3)$ or $p_3(\cdot,I,\cdot)$. There are a total\footnote{Resulting from all the possible curryings of 4 identity cubices in each of the 3 arguments.} of 12 different such binary operations in each cubix algebra. To make matters more interesting, most of them satisfy no obvious axioms, in particular they are non-associative, such as $p_3(\cdot,I,\cdot)$, while some of them turn out to be very well-behaved, as is the case for $p_1(\cdot,\cdot,i_3)$ which is associative and unital.\newline

Lastly, we turn to the question of generalized associativity: what is the analogue of binary matrix multiplication associativity for ternary multiplications of cubices? More specifically, we are interested in knowing whether compositions of multiple ternary multiplications satisfy some generic properties such as
\begin{align*}
    & p_1(p_1(a,b,c),d,e) \stackrel{?}{=} p_1(a,p_1(b,c,d),e) \stackrel{?}{=} p_1(a,b,p_1(c,d,e))\\
    & p_2(a,p_2(b,c,d),p_2(e,f,g)) \stackrel{?}{=} p_2(p_2(a,b,d),c,p_2(e,f,g))\\
    & \text{etc.}
\end{align*}
This question was partially addressed in \cite{abramov2009algebras,abramov2011ternary}, whence it was found that no ternary multiplication of cubices satisfies na\"ive ternary associativity (first line of equalities above) and that $p_4$ satisfies the so-called \emph{ternary associativity of the second kind}:
\begin{equation*}
    p_4(p_4(a,b,c),d,e) = p_4(a,p_4(d,c,b),e) = p_4(a,b,p_4(c,d,e)).
\end{equation*}
This axiom corresponds to the associativity-like property of heaps \cite{kolar2000heap,santiago2010ternary}, which are ternary algebraic structures that have been known since the early 20$^\text{th}$ century. We suspect that the relatively well-developed understanding of the fish algebra found in the existing literature is due to the strongly sequential nature of the fish multiplication. Ongoing work by the authors \cite{zapata2022heap} indeed suggests that there is an elegant argument to relate the heap axiom to the hypergraphic properties of the fish shape. Using the computer algebra implementation we were able to check that $p_1$, $p_2$, $p_3$ satisfy no axioms involving 5 or 7 elements. Axioms involving 9 or more elements were not checked due to the limited computing power available to the authors. For all we know, the $p_1$, $p_2$ and $p_3$ multiplications may satisfy relatively intelligible associativity-like properties involving 9 elements, they may satisfy obscure complicated properties involving large amounts of elements or they may simply not be amenable to this kind of axiom and fail to satisfy any form of higher associativity. The mystery remains.\newline

In our view, the question of ternary associativity perfectly encapsulates the current impasse in the understanding of higher arity structures. We seem to have a confident understanding of the ingredients involved, after all, cubices are just arrays of numbers and their multiplications are defined by simple arithmetical operations, and we have grasped the ternary generalization of concepts such as composition and identity. Yet, something as apparently elementary as the analogue of associativity defies our intuition and eludes a concrete definition. Making any significant progress on this front will surely result in a qualitative leap towards the general understanding of higher arity.

\subsection{Ternary Relations} \label{rel}

Higher order relations have been occasionally studied in the context of logic and set theory since the mid 19$^\text{th}$ century, as seen in the pioneering work of C. S. Peirce \cite{peirce1870description,peirce1902logic} containing explicitly higher arity ideas. However,  developments in logic \cite{peirce1880algebra}, the theory of relations \cite{tarski1941calculus,fraisse2000theory} and formal languages \cite{moll2012introduction} focused almost exclusively on binary relations and sequential constructions. We do find select research pieces on relational structures \cite{mcculloch1967triadas,novak1989transitive,konecny2014triadic,zedam2018traces,bakri2021compositions}, association schemes \cite{mesner1990association} and semantics \cite{beall2012ternary} that have investigated concrete higher arity relations in detail, although most work in these areas is often directed towards generalist questions such as algorithmic complexity \cite{kumar1992algorithms} or universality \cite{topentcharov1993composition,behrisch2013relational}. Our goal is to focus on binary and ternary relations to highlight their elementary compositional properties.\newline

A binary relation between two sets $A$ and $B$ is a set of pairs $R\subset A\times B$. Diagrammatically, we represent a relation as an edge in a graph whose vertices represent the related sets. When elements $a\in A$ and $b\in B$ are related we write $R_{ab}$, that is, $ (a,b)\in R \Leftrightarrow R_{ab}$. Two relations can be composed if they share a common set: let two relations $R\subset A\times B$ and $S\subset B\times C$, the composed relation is defined by middle-element transitivity
\begin{equation*}
    S \circ R :=\{(a,c)\in A\times C\,| \, \exists\, b\in B,\, R_{ab}\, \wedge \, S_{bc}\}.
\end{equation*}
Composition of binary relations is associative by construction and the diagonal relations $i_A:=\{(a,a))\,|\, a\in A \}$ act as compositional identities for any set $A$. A \textbf{ternary relation} is simply defined as a set of trios of elements from three sets $T\subset A\times B \times C$ and we shall write $ (a,b,c)\in T \Leftrightarrow T_{abc}$. Ternary relations are thus represented as triangular hyperedges in the relational graph. In general, $n$-ary relations are defined as sets of $n$-tuples and are represented as hyperedges of arity $n$. Note that our order agnosticism clashes here with the need to use cartesian products of sets to define relations. We shall again leave out any matters of order and sequentiality for ternary relations.\newline 

It is easy to see that compositions of $n$-ary relations via middle-element transitivity, generalizing the binary composition formula above, and middle-element equality are formally indistinguishable from index contractions and incidences of $n$-matrices as defined in \ref{cubix}. It follows from this observation that the compositional theory of ternary relations can be developed following the template provided by the compositional framework of $3$-graphs. The discussion in \ref{cubix} applies here mutatis mutandis replacing $\mathbb{F}^N$ spaces by general sets and $3$-matrices by ternary relations.\newline

Following Figure \ref{graph3}, \textbf{ternary arity-preserving compositions} are defined analogously to the $3$-matrix multiplications $p_1$--$p_7$ in \ref{cubix}. We give here the explicit set-theoretic formulas for the compositions induced by the cone, blades, triforce and fish shapes:
\begin{align*}
    \triangle_1 (R,S,T) &:=\{(a,b,c)\in A\times B\times C\,| \, \exists\, x\in X,\ R_{abx}\, \wedge \, S_{axc}\, \wedge \, T_{xbc}\}  &\text{(cone)}\\
    \triangle_2 (R,S,T) &:=\{(a,b,c)\in A\times B\times C\,| \, \exists\, x\in X,\, y\in Y,\, R_{axy}\, \wedge \, S_{xby}\, \wedge \, T_{xyc}\} &\text{(blades)}\\
    \triangle_3 (R,S,T) &:=\{(a,b,c)\in A\times B\times C\,| \, \exists\, x\in X,\, y\in Y,\, z\in Z,\, R_{axy}\, \wedge \, S_{xbz}\, \wedge \, T_{xzc}\} &\text{(triforce)}\\
    \triangle_4 (R,S,T) &:=\{(a,b,c)\in A\times B\times C\,| \, \exists\, x\in X,\, y\in Y,\, z\in Z,\, R_{abx}\, \wedge \, S_{xyz}\, \wedge \, T_{yzc}\} &\text{(fish)}
\end{align*}
where $R,S,T$ are arbitrary ternary relations among the sets $A$, $B$, $C$, placed in external vertices in the hypergraph diagrams, and sets $X$, $Y$, $Z$, placed in internal vertices. Composition $\triangle_1$ was defined indirectly in \cite{mesner1990association} as part of a ternary association scheme and $\triangle_3$, $\triangle_4$ were defined in \cite{mcculloch1967triadas} in a pioneering attempt to understand compositions of ternary relations. To our knowledge, $\triangle_2$, $\triangle_5$, $\triangle_6$, $\triangle_7$ are novel and introduced here for the first time.\newline

The question of ternary associativity remains as mysterious for ternary relations as it appeared for $3$-matrices in \ref{cubix}. Nevertheless, it turns out that ternary relations display the properties of ternary identities much more transparently than $3$-matrices. Ternary identity relations are defined by partial diagonal subsets $i_{AAB}:=\{(a,a,b)\,|\, a\in A \,, b\in B\}\subset A\times A \times B$, called \textbf{partial identity relations}, and diagonal subsets $I_A:=\{(a,a,a)\,|\, a\in A\}\subset A\times A\times A$, called \textbf{tridentity relations}. Not only are these definitions direct generalizations of the binary notion of identity relation as diagonal subset but they also operate diagrammatically in entirely analogous ways. To see this, note that we can easily adapt the matrix splicing notation introduced in \ref{cubix} to represent relations between sets: a binary relation $R\subset A\times B$ and a ternary relation $T\subset A\times B \times C$ are simply notated as
\begin{center}
    \includegraphics[scale=0.17]{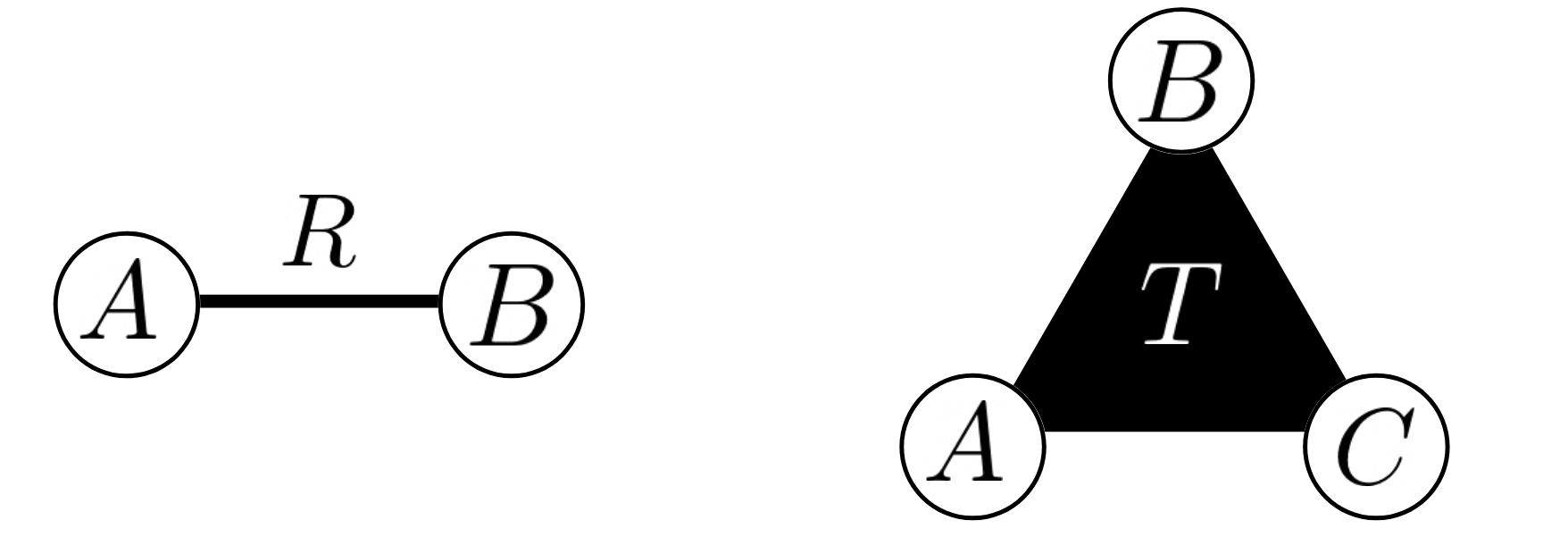}
\end{center}
The compositional behaviour of binary identity relations can be efficiently captured in the following diagrammatic property
\begin{center}
    \includegraphics[scale=0.19]{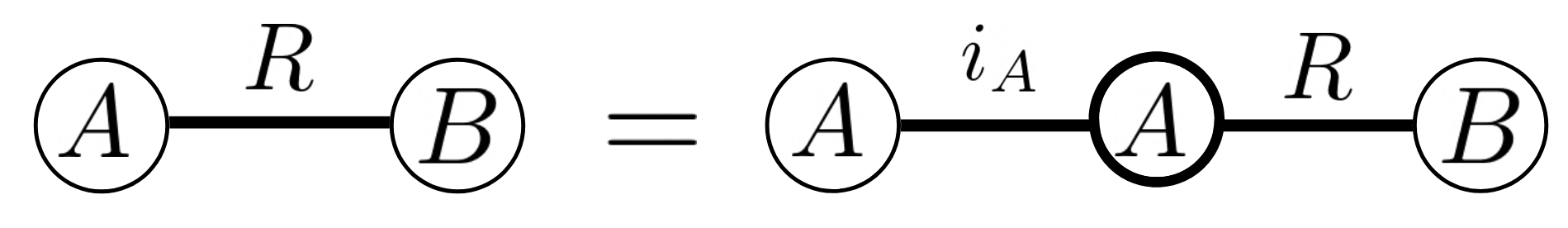}
\end{center}
Similarly, two partial identities simplify a cone composition:
\begin{center}
    \includegraphics[scale=0.19]{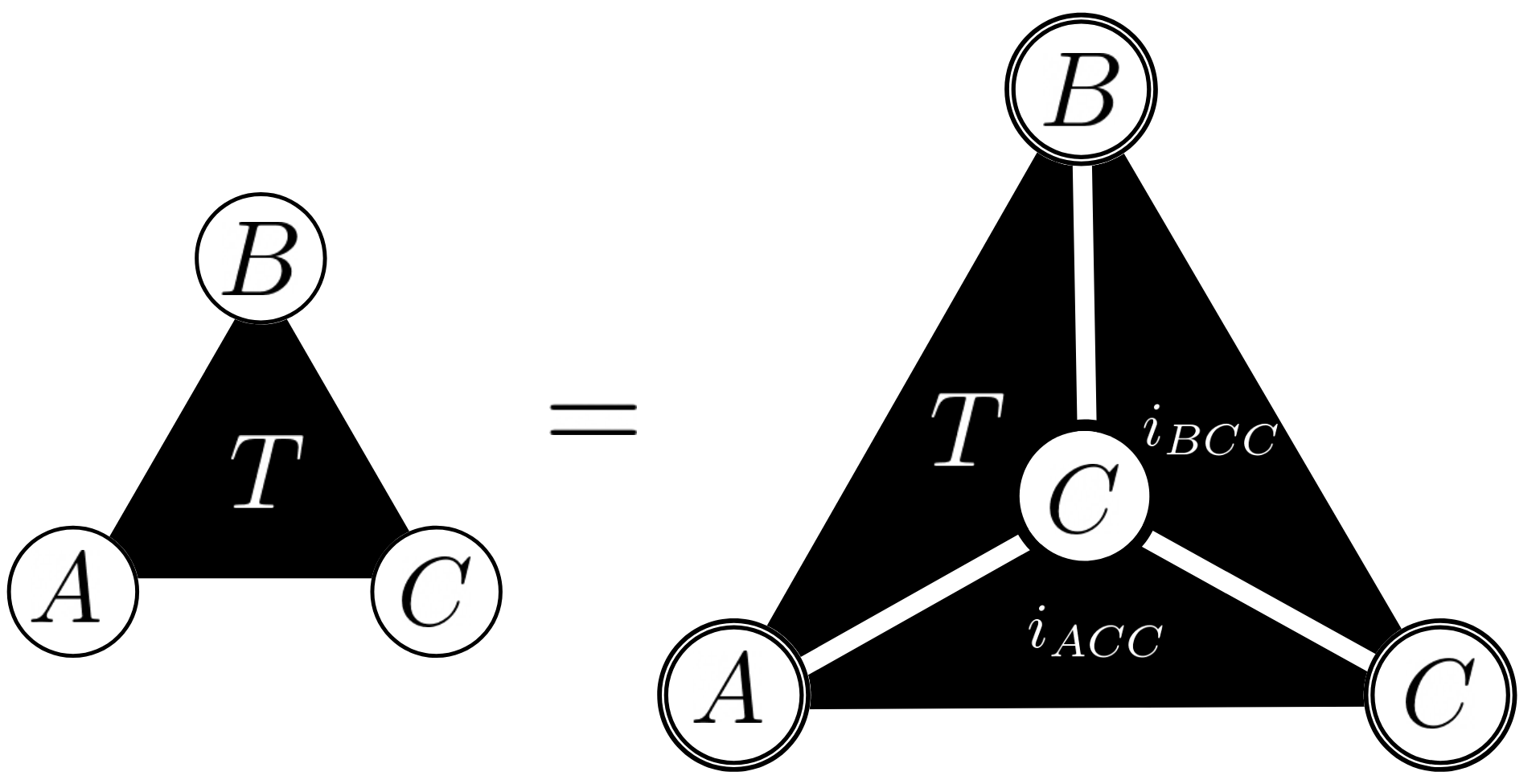}
\end{center}
and a partial identity together with a tridentity simplify a triforce composition:
\begin{center}
    \includegraphics[scale=0.19]{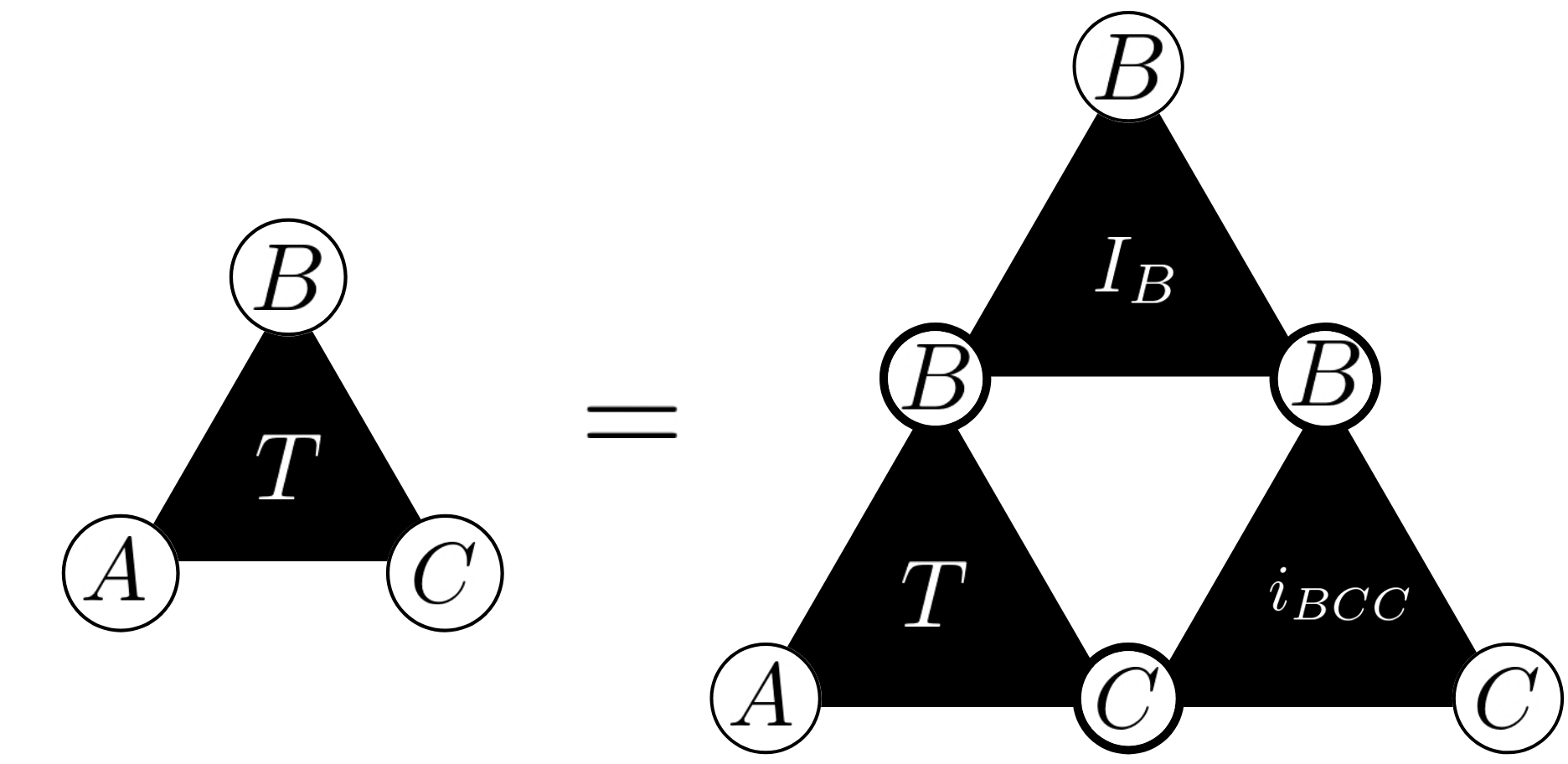}
\end{center}

\subsection{Ternary Categories} \label{cat}

The basic definition of a category \cite{mac2013categories,leinster2014basic} has a strongly binary and sequential flavour: morphisms are directed assignments between pairs of objects. Throughout the second half of the 20$^\text{th}$ century, several developments in category theory gradually introduced higher notions that extended the basic binary structure of categories: multicategories \cite{lambek1989multicategories}, $n$-categories \cite{bastiani1974multiple}, weak $n$-categories \cite{benabou1967introduction}, polycategories \cite{szabo1975polycategories}, operads \cite{kriz1995operads} and opetopes \cite{baez1998higher}. Despite the manifestly higher order nature of these generalizations, they all carry the intrinsic binary nature of arrow-like morphisms at their core. The recent work of N. Baas \cite{baas2019mathematics,baas2019philosophy} underlines this aspect of typical categorical formalisms and addresses its limitations by proposing so-called hyperstructures, which serve as a general template that will fit higher arity generalizations of categories. The first (and, to our knowledge, only) attempt to systematically investigate category-like structures whose morphisms relate general collections of objects can be found in the work of V. Topentcharov in the late 20th century \cite{topentcharov1988structures,topentcharov1993composition}. The concluding remarks of \cite{topentcharov1988structures} include a mention of $3$-matrix algebras and ternary relations as putative examples of ternary categories to be studied in detail in subsequent work. To this day, it seems that there has been no further research picking up that thread.\newline

Similar to how the class of binary relations on sets serves as a model for the definition of ordinary categories \cite{leinster2014basic}, our discussion on ternary relations in \ref{rel} offers a glimpse at what ternary categories may look like. Given the vast scope and high levels of abstraction of modern algebra and category theory \cite{wiki2008nlab}, it may seem surprising at first that there is no agreed-upon definition of ternary category. Although we suspect that the binarity bias of \ref{binary} is at play again here as a broad cause of the status quo, given our investigation into ternary operations in \ref{cubix} and \ref{rel}, we can nevertheless confidently attribute the lack of a definition of ternary category to the poor understanding of ternary associativity. As soon as associativity-like properties are successfully captured for ternary objects such as $3$-matrices or ternary relations, natural definitions of ternary category are likely to follow easily.\newline

Despite the fundamentally binary nature of ordinary categories, we can exploit their basic properties to construct ternary structures. We shall see that this is not dissimilar to how in \ref{cubix} and \ref{rel} we defined genuinely ternary operations of matrices and relations despite relying on strongly binary and sequential structures such as number fields, arrays, cartesian products, Boolean algebra, etc. Our goal is to define the ternary analogue of the elementary notion of isomorphism in a category $\mathcal{C}$. Recall that an isomorphism between a pair of objects $A, B \in \mathcal{C}$ is defined as a pair of morphisms fitting in the following commutative diagram, where commutativity amounts to the property of inverse:
\begin{equation*}
\begin{tikzcd}[]
A \arrow[rrr, bend left, "f"] & & & B \arrow[lll, bend left, "g"]
\end{tikzcd}
\qquad \qquad
\begin{array}{cc}
g \circ f = 1_A\\
f \circ g = 1_B
\end{array}
\end{equation*}
We can easily extend this notion to a finite family of $n$ objects and define an \textbf{$n$-ary isomorphism} as a commuting cycle of morphisms between them. In particular, due to our interest in the ternary case we define a \textbf{ternary isomorphism}, or \textbf{trisomorphism} for short, between a trio of objects $A, B, C \in \mathcal{C}$ as a trio of morphisms fitting in the following commutative diagram:
\begin{equation} \label{triso}
\begin{tikzcd}[baseline=(current bounding box.center)]
                             & B \arrow[dr, bend left, "g"] &  \\
A \arrow[ur, bend left, "f"] &                              & C \arrow[ll, bend left, "h"]
\end{tikzcd}
\qquad \qquad
\begin{array}{cc}
h \circ g \circ f = 1_A \\
f \circ h \circ g = 1_B \\
g \circ f \circ h = 1_C
\end{array}
\tag{triso}
\end{equation}
It follows from this definition that each of the morphisms in a trisomorphism is also a binary isomorphism, however not any triangle of binary isomorphisms between a trio of objects forms a trisomorphism. Although trisomorphisms are determined by any pair of the trio of morphisms we will notate all three explicitly to highlight the ternary nature of the constructions we propose below. In this sense, we shall regard trisomorphisms as a manner of ternary morphisms between objects.\newline

Binary isomorphisms in ordinary categories can be represented by graph-like diagrams, as directedness of morphisms is lost due to the existence of inverses. Similarly, we may represent trisomorphisms by $3$-graphs where the vertices are objects and the hyperedges are given by commuting cycles of trios of morphisms. We are looking for elementary ways to define trisomorphism compositions. Let us begin by considering the simple case of two trisomorphisms sharing two objects:
\begin{equation*}
    \begin{tikzcd}[row sep=small, column sep=large]
 & & X \arrow[dd, bend right, "g_1"'] \arrow[drr, bend left, "h_2"] & & \\
A \arrow[urr, bend left, "f_1"] & & & & B \arrow[dll, bend left, "f_2"] \\
 & & Y \arrow[uu, bend right, "g_2"'] \arrow[ull, bend left, "h_1"] & & 
    \end{tikzcd}
\end{equation*}
This diagram suggests that two trisomorphisms arranged this way may be spliced somehow to define a binary isomorphism between $A$ and $B$. Indeed, the morphisms $\alpha=h_2\circ f_1$ and $\beta=h_1 \circ f_2$ map between the correct objects but the trisomorphism condition alone does not imply that $\alpha$ and $\beta$ are inverses. A sufficient condition for this to happen is the natural requirement of \textbf{isomorphism uniqueness}: in any arrangement of a collection of $n$-ary isomorphisms, whenever more than one binary isomorphisms occur between a pair of objects they must all be equal. In the case at hand this amounts to $g_1=g_2^{-1}$ and it is a simple check to verify that this, together with the trisomorphism condition, implies that $\alpha=\beta^{-1}$.\newline

Given a collection of higher arity isomorphisms sharing some objects, we define a general \textbf{splicing} as a new higher arity isomorphism constructed from compositions of individual morphisms when the isomorphism uniqueness condition holds between pairs of shared objects. We can see, once more, that plausible splicings of higher arity isomorphisms are suggested by hypergraph splicings defined \ref{graph}. Concretely for trisomorphisms, Figure \ref{graph3} offers the template to define the simplest splicings. By parsing $3$-graph shapes in different ways we find many possibilities to construct ternary \textbf{compositions} of trisomorphisms, i.e. arity-preserving splicings involving three trisomorphisms. Given their particular ternary symmetry and the attention we devoted to them in previous sections, we shall focus on three particular constructions using the cone, blades and triforce shapes.\newline

The definition of ternary compositions of trisomorphisms can be easily understood by a straightforward diagrammatic procedure: we begin with three objects $A$, $B$, $C$ that we label as outer, then progressively throw in inner objects $X$, $Y$, $Z$ linked by trisomorphisms with the outer objects in a pattern so as to form the cone, blades and triforce shapes. The goal is to get a trisomorphism diagram between the outer objects as in (\ref{triso}) by imposing the necessary isomorphism uniqueness conditions. We shall draw compositional diagrams following simplicial-like orientations of labelled morphisms for ease of notation. The \textbf{cone composition} is given by the diagram
\begin{equation*}
\begin{tikzcd}[row sep=large, column sep=large]
& & B \arrow[ddrr, bend left, "f_2"] \arrow[d, bend right, "g_1"'] & & \\[10pt]
& & X \arrow[dll, bend right, "h_1"'] \arrow[u, bend right, "h_2"'] \arrow[drr, bend right, "h_3"'] & & \\[-10pt]
A \arrow[uurr, bend left, "f_1"] \arrow[urr, bend right, "g_3"'] & & & & C \arrow[llll, bend left, "f_3"] \arrow[ull, bend right, "g_2"']
\end{tikzcd}
\end{equation*}
Here, isomorphism uniqueness demands that the morphisms from the outer objects to $X$ are all the same, that is $g_1=h_2^{-1}$, $g_2=h_3^{-1}$, $g_3=h_1^{-1}$, then the outer trisomorphism is simply given by $\alpha=f_1$, $\beta=f_2$, $\gamma=f_3$. The \textbf{blades composition} results from the diagram
\begin{equation*}
\begin{tikzcd}[]
 & &[-4pt] B \arrow[dl, bend right, "f_2"'] & & \\[10pt]
 & X \arrow[from=drrr, bend left, near start, "f_3"'] \arrow[rr, bend left, "g_1"] \arrow[rr, "g_2" description] \arrow[rr, bend right,"g_3"'] & & Y \arrow[dlll, bend left, crossing over, near end, "h_1"'] \arrow[dr, bend left, "h_3"] \arrow[ul, bend right, "h_2"'] & \\[6pt]
A \arrow[ur, bend left, "f_1"] & & & & C
\end{tikzcd}
\end{equation*}
Here, isomorphism uniqueness demands that the morphisms between $X$ and $Y$ are all equal $g_1=g_2=g_3=:z$ and the outer trisomorphism is given by $\alpha=h_2 \circ z \circ f_1$, $\beta= h_3 \circ z \circ f_2$, $\gamma=h_1 \circ z \circ f_3$. The \textbf{triforce composition} is given by the diagram
\begin{equation*}
\begin{tikzcd}[]
 & & B \arrow[dr, bend left, "f_2"] & & \\[5pt]
 & X \arrow[dr, bend right, "g_1"] \arrow[ur, bend left, "h_2"] & & Z \arrow[dr, bend left, "h_3"] \arrow[ll, bend right, "g_2"] & \\
A \arrow[ur, bend left, "f_1"] & & Y \arrow[ll, bend left, "h_1"] \arrow[ur, bend right, "g_2"] & & C \arrow[ll, bend left, "f_3"]
\end{tikzcd}
\end{equation*}
Here, isomorphism uniqueness demands that multiple paths between outer objects compose to the same morphism, which turns out to be equivalent to the condition that the triangle of morphisms between the inner objects is a trisomorphism. Explicitly, the composability condition is given by $g_1=g_2^{-1} \circ g_3^{-1}$, $g_2=g_3^{-1} \circ g_1^{-1}$, $g_3=g_1^{-1} \circ g_2^{-1}$, so that the outer trisomorphism is then defined by $\alpha = h_2\circ f_1$, $\beta = h_3 \circ f_2$, $\gamma = h_1 \circ f_3$.\newline

The existence of identity morphisms on all objects of a category induces notions of ternary identities similar to those described for compositions of ternary relations in \ref{rel}. Any binary isomorphism induces a \textbf{partial identity trisomorphism} by trivially inserting an identity arrow to complete the triangle:
\begin{equation*}
\begin{tikzcd}[]
A \arrow[rrr, bend left] \arrow[rrr, phantom, "f" description] & &[20pt] & B \arrow[lll, bend left]
\end{tikzcd}
\,
\text{ isomorphism }
\quad \Leftrightarrow \quad
\begin{tikzcd}[baseline=(current bounding box.center)]
                             & A \arrow[dr, bend left, "f"] &  \\
A \arrow[ur, bend left, "1_A"] &                              & B \arrow[ll, bend left, "f^{-1}"]
\end{tikzcd}
\,
\text{ trisomorphism }
\end{equation*}
Similarly, the \textbf{tridentity trisomorphism} can be constructed on any object:
\begin{equation*}
\begin{tikzcd}[baseline=(current bounding box.center)]
                             & A \arrow[dr, bend left, "1_A"] &  \\
A \arrow[ur, bend left, "1_A"] &                              & A \arrow[ll, bend left, "1_A"]
\end{tikzcd}
\end{equation*}
In order to clearly distinguish partial identities and tridentities in categorical diagrams we use equal signs:
\begin{equation*}
\text{ partial identity: } \,
\begin{tikzcd}[baseline=(current bounding box.center)]
                             & A \arrow[dr, bend left, "f"] &  \\
A \arrow[ur, bend left, equals] &                              & B \arrow[ll, bend left, "f^{-1}"]
\end{tikzcd}
\qquad \qquad
\text{ tridentity: } \,
\begin{tikzcd}[baseline=(current bounding box.center)]
                             & A \arrow[dr, bend left, equals] &  \\
A \arrow[ur, bend left, equals] &                              & A \arrow[ll, bend left, equals]
\end{tikzcd}
\end{equation*}
The compositional behaviour of these special trisomorphisms indeed appears identity-like in a diagrammatic sense when considering the compositions defined above. The cone composition satisfies a compatibility condition between tridentity and partial identities
\begin{equation*}
\begin{tikzcd}[]
A \arrow[rrr, bend left] \arrow[rrr, phantom, "f" description] & &[20pt] & B \arrow[lll, bend left]
\end{tikzcd}
\quad = \quad
\begin{tikzcd}[row sep=large, column sep=large,baseline=(current bounding box.center)]
& & B \arrow[ddrr, bend left, equals] & & \\[10pt]
& & B \arrow[dll, bend right, near start, "f^{-1}"'] \arrow[u,equals] \arrow[drr, equals] & & \\[-10pt]
A \arrow[uurr, bend left, "f"] \arrow[urr, bend right, "f"'] & & & & B \arrow[llll, bend left, "f^{-1}"] 
\end{tikzcd}
\end{equation*}
and a general simplifying property for partial identities:
\begin{equation*}
\begin{tikzcd}[baseline=(current bounding box.center)]
                             & B \arrow[dr, bend left, "g"] &  \\
A \arrow[ur, bend left, "f"] &                              & C \arrow[ll, bend left, "h"]
\end{tikzcd}
\quad = \quad
\begin{tikzcd}[row sep=large, column sep=large,baseline=(current bounding box.center)]
& & B \arrow[ddrr, bend left, "g"] & & \\[10pt]
& & B \arrow[dll, bend right, near start, "f^{-1}"'] \arrow[u,equals] \arrow[drr, bend right, "g"'] & & \\[-10pt]
A \arrow[uurr, bend left, "f"] \arrow[urr, bend right, "f"'] & & & & C \arrow[ull, bend right, near end, "g^{-1}"'] \arrow[llll, bend left, "h"] 
\end{tikzcd}
\end{equation*}
Similarly, an identities compatibility property holds for the blades composition
\begin{equation*}
\begin{tikzcd}[]
A \arrow[rrr, bend left] \arrow[rrr, phantom, "f" description] & &[20pt] & B \arrow[rrr, bend left] \arrow[lll, bend left] \arrow[rrr, phantom, "g" description] & &[20pt] & C \arrow[lll, bend left]
\end{tikzcd}
\quad = \quad
\begin{tikzcd}[]
 & &[-4pt] B \arrow[dl, bend right, equals] & & \\[10pt]
 & B \arrow[from=drrr, bend left, near start, "g^{-1}"'] \arrow[rr, "g" description] & & C \arrow[dlll, bend left, crossing over, near end, "(g\circ f)^{-1}"] \arrow[dr, bend left, equals] \arrow[ul, bend right, "g^{-1}"'] & \\[6pt]
A \arrow[ur, bend left, "f"] & & & & C
\end{tikzcd}
\end{equation*}
and the triforce composition
\begin{equation*}
\begin{tikzcd}[row sep=large, column sep=large]
& & B \arrow[d, bend right] \arrow[d, phantom, "g" description] & & \\[10pt]
& & X \arrow[dll, bend right] \arrow[u, bend right] \arrow[drr, bend right] & & \\[-10pt]
A \arrow[urr, bend right] \arrow[urr, phantom, "f"] & & & & C \arrow[ull, bend right] \arrow[ull, phantom, "h" description]
\end{tikzcd}
\quad = \quad
\begin{tikzcd}[]
 & & X \arrow[dr, bend left, equals] & & \\[5pt]
 & X \arrow[dr, bend right, "g^{-1}"] \arrow[ur, bend left, equals] & & X \arrow[dr, bend left, "h^{-1}"] \arrow[ll, bend right, equals] & \\
A \arrow[ur, bend left, "f"] & & B \arrow[ll, bend left, "f^{-1} \circ g"] \arrow[ur, bend right, "g"] & & C \arrow[ll, bend left, "g^{-1} \circ h"]
\end{tikzcd}
\end{equation*}
Lastly, a general simplifying property also holds for the blades composition
\begin{equation*}
\begin{tikzcd}[baseline=(current bounding box.center)]
                             & B \arrow[dr, bend left, "g"] &  \\
A \arrow[ur, bend left, "f"] &                              & C \arrow[ll, bend left, "h"]
\end{tikzcd}
\quad = \quad
\begin{tikzcd}[]
 & &[-4pt] B \arrow[dl, bend right, equals] & & \\[10pt]
 & B \arrow[from=drrr, bend left, near start, "g^{-1}"'] \arrow[rr, "g" description] & & C \arrow[dlll, bend left, crossing over, near end, "h"'] \arrow[dr, bend left, equals] \arrow[ul, bend right, "g^{-1}"'] & \\[6pt]
A \arrow[ur, bend left, "f"] & & & & C
\end{tikzcd}
\end{equation*}
and the triforce composition
\begin{equation*}
\begin{tikzcd}[baseline=(current bounding box.center)]
                             & B \arrow[dr, bend left, "g"] &  \\
A \arrow[ur, bend left, "f"] &                              & C \arrow[ll, bend left, "h"]
\end{tikzcd}
\quad = \quad
\begin{tikzcd}[baseline=(current bounding box.center)]
 & & B \arrow[dr, bend left, equals] & & \\[5pt]
 & B \arrow[dr, bend right, "g"] \arrow[ur, bend left, equals] & & B \arrow[dr, bend left, "g"] \arrow[ll, bend right, equals] & \\
A \arrow[ur, bend left, "f"] & & C \arrow[ll, bend left, "h"] \arrow[ur, bend right, "g^{-1}"] & & C \arrow[ll, bend left, equals]
\end{tikzcd}
\end{equation*}
Note the similarities between these diagrammatic equalities and the compositional properties of identity ternary relations shown at the end of \ref{rel}.\newline

Although so far in our research on trisomorphisms we have not been able to make any further progress on the question of ternary associativity, being able to work with precisely defined notions in categories poses a great advantage over other set-theoretic approaches, such as cubix algebras \ref{cubix} and ternary relations \ref{rel}, which tend to be formally heavier and slightly more ambiguously defined. We currently place our hope on the fact that trisomorphism compositions are ultimately defined in terms of binary morphisms, and so one expects binary associativity to somehow permeate upwards into ternary and higher arity compositions.

\subsection{$3$-Lie Algebras} \label{Lie}

Among all the mathematical structures that generalize binarity, higher arity algebras are perhaps the ones that have received the most attention \cite{lister1971ternary,michor1998n,carlsson1980n,markl2002operads,arnlind2010ternary,gal2011n,cerejeiras2021ternary}. The sustained interest in $n$-ary algebras is likely to have been fueled by the developments in theoretical and mathematical physics discussed in \ref{phys}. This has been particularly true for $n$-ary Lie algebras, also known as Filippov algebras \cite{filippov1985n}, as explained in the excellent review \cite{azcarraga2010nary}. We shall focus on the ternary instance of these algebras: a \textbf{$3$-Lie algebra} is a vector space endowed with a trilinear operation $(\mathfrak{g},[\,\,,\,,\,])$ that is totally skew-symmetric, i.e. $[\,\,,\,,\,] : \wedge^3 \mathfrak{g} \to \mathfrak{g}$, and that satisfies the \textbf{ternary Jacobi identity}, sometimes also called the fundamental identity:
\begin{equation*}
    [x,y,[a,b,c]]=[[x,y,a],b,c]+[a,[x,y,b],c]+[a,b,[x,y,c]].
\end{equation*}

Ordinary ($2$-)Lie algebras play a central role in modern physical theories since they encode the infinitesimal symmetry of Lie groups \cite{azcarraga1998lie,hall2003lie}. A celebrated result in the theory of Lie groups and Lie algebras is the so-called Lie's Third Theorem \cite{cartan1952theorie} which states that for any finite-dimensional Lie algebra there exists a unique simply-connected Lie group that integrates it. This means that, aside from global topological aspects, Lie algebras and Lie groups are conceptually equivalent. A natural question then arises: \textbf{is there an analogue of Lie's Third Theorem for $3$-Lie algebras?} As it stands, this is an open problem \cite{azcarraga2010nary}.\newline

That this question remains unresolved to this day is perplexing for several reasons. Firstly, given the highly sophisticated mathematical theories and the immense computational power at our disposal, one may think that the relatively simple axioms of a $3$-Lie algebra would pose no significant challenge. But perhaps more puzzling is the fact that, among the vast diversity of group-like algebraic structures that have been defined over the years, one finds no apparent candidates for $3$-Lie algebra integrators. Higher arity groups with na\"ively generalized associativity axioms \cite{post1940polyadic,dudek2007remarks,gal2011n}, hypergroups with operations defined on subsets instead of single elements \cite{davvaz2009generalization,massouros2021overview} and the $n$-categorical generalizations of groups \cite{brown1981equivalence,zhu2006lie}, despite some recent partial results \cite{huerta2015division}, all fail to act as $3$-Lie algebra integrators in any significant degree of generality.\newline

What makes the problem of $3$-Lie algebra integration particularly enticing is that simple examples of $3$-Lie algebras abound. The most natural instance of $3$-Lie algebra occurs in $4$-dimensional euclidean space: let $\{\mathbf{e}_i\}_{i=1}^4$ be a basis in $\Real^4$ (regarded as an euclidean vector space with the standard inner product) and take three generic elements $a,b,c\in\Real^4$, then the manifestly skew-symmetric ternary bracket defined by
\begin{equation*}
[a,b,c]:=
    \begin{vmatrix}
    \mathbf{e}_1 & \mathbf{e}_2 & \mathbf{e}_3 & \mathbf{e}_4 \\
    a_1 & a_2 & a_3 & a_4 \\
    b_1 & b_2 & b_3 & b_4 \\
    c_1 & c_2 & c_3 & c_4
    \end{vmatrix} 
\end{equation*}
can be easily checked to satisfy the ternary Jacobi identity. Note that this is the canonical codimiension-$1$ orthogonal complement operation present in all euclidean spaces: in $\Real^2$ it manifests as the complex structure or rotation by $\pi/2$, in $\Real^3$ as the usual cross product Lie algebra $\mathfrak{so}(3)$ and in $\Real^{n+1}$ as a $n$-Lie algebra structure in general. These so-called \textbf{euclidean $n$-Lie algebras} have received considerable attention in the mathematical physics literature due to the famous BLG model \cite{bagger2008gauge} that used a $3$-Lie bracket to describe M2 brane interactions. Another example of a $3$-Lie algebra was found in efforts to quantize Nambu-Poisson structures \cite{awata2001quantization}: for three ordinary square matrices $A,B,C\in \mathbb{F}^{N^2}$ the bracket
\begin{equation*}
    [A,B,C]:=\mathrm{tr}(A)[B,C]+\mathrm{tr}(B)[C,A]+\mathrm{tr}(C)[A,B]
\end{equation*}
is skew-symmetric and satisfies the ternary Jacobi identity. Similar constructions of $3$-Lie brackets have also been found for cubices \cite{bai20143lie}. $n$-Lie algebras have been related to their higher-categorical analogous Lie $n$-algebras \cite{chen2021lie}. Despite the dangerously similar notation, Lie $n$-algebras are a particular kind of L$_\infty$-algebras, which are structures generalizing ordinary ($2$-)Lie algebras via higher homotopy encompassing many important objects in differential geometry and mathematical physics \cite{stasheff2019infinity}.\newline

A first approach to the general problem of $3$-Lie algebra integration is to simply search for candidates of integrators: can we find a smooth manifold carrying a ternary operation, presumably with appropriate notions of identity and invariance, which somehow induces an infinitesimal $3$-Lie algebra structure on (some submodule of sections of) the tangent bundle? By analogy with the case of $2$-Lie algebra integration, we would expect the integrating ternary structure to be group-like, satisfying some form of ternary associativity. If this turns out to be the case, the lack of integrator candidates comes as no surprise after our discussion on ternary associativity in previous sections. We thus conclude that a better understanding of higher associativity is likely to result in substantial progress towards the long-standing problem of $n$-Lie algebra integration.\newline 

The difficulties resulting from the lack of general understanding of ternary structures are compounded by a subtle issue already present in the binary case: J.L. Loday's \emph{coquecigrue} problem \cite{loday1993version}. Leibniz algebras are non-commutative analogues of Lie algebras, that is, vector spaces with a bilinear operation that satisfies the Jacobi identity but has no particular symmetry properties. Leibniz algebras have been shown to be integrated by Lie racks \cite{monterde2014integral}, which are smooth manifolds carrying a compatible binary operation that satisfies a distributivity property instead of associativity. All Lie groups have an associated Lie rack given by the conjugation operation. The \emph{coquecigrue} problem emerges because one obtains two inequivalent results when integrating a Lie algebra with a Lie group (in the usual fashion) and with a Lie rack (regarding the Lie algebra as a Leibniz algebra). This is particularly relevant to the ternary case because an elementary way to obtain a ternary operation is by nesting binary operations. In the case of a Lie algebra $(\mathfrak{g},[\,,])$, the ternary nested bracket $[ \,,[\,,]]$ satisfies the ternary Jacobi identity but it is not totally skew-symmetric. This is an instance of a \textbf{$3$-Leibniz algebra}, more concretely a Lie triple system \cite{lister1952structure}. As in the binary case, a complete solution to the $3$-Lie algebra integration problem would also need to account for $3$-Leibniz algebras.\newline

Recent work on non-associative ternary algebras \cite{kerner2008ternary,abramov2011ternary,abramov2017algebras} shows that the symmetry properties of ternary operations should not be limited to the usual skew-symmetrization based on the $\pm 1$ signature of permutations, and invites to consider more general representations of $\mathbb{Z}_3$ and $S_3$. This suggests that the definition of $3$-Lie algebras with totally skew-symmetric brackets is likely to represent only a small family within the zoology of $3$-Leibniz algebras. Beyond symmetry considerations, there is increasing evidence that points to other ternary generalizations of the Jacobi identity. For instance, motivated by the ternary commutator of bounded operators on a Hilbert space, a modified version of the ternary derivation property appeared in a recent paper examining the role ternary algebras in quantum mechanics \cite{bruce2021semiheaps}.

\subsection{Topology and Geometry} \label{topo}

The Borromean rings were presented in \ref{binary} as the most compact and intuitive manifestation of irreducible ternarity. Beyond the knot theoretic generalization of the Borromean rings, known as Brunnian links \cite{debrunner1961links,miyazawa2006classification}, the fields of topology and geometry are rich sources of higher arity. We suspect that the spatial nature of these branches of mathematics somehow keeps binarity at bay and favours the naturally higher order notions of shape, locus and structure. A synthesis of these ideas is beautifully captured in the work of N. Baas on higher topological structures \cite{baas2016higher}.\newline

Although we postpone discussing any concrete ternary topological structures for future work, we would like to remark that topology offers a highly developed toolkit of concepts and mathematical theories that may aid in the systematic investigation of higher arity. Particularly, the theory of simplicial complexes and topological manifolds \cite{lee2010introduction} provides some useful notions to begin a systematic enquiry into ternary sequentiality. A sequence is a very special kind of simplicial complex, one whose pairs of $0$-simplices are all connected by a single $1$-simplex. In other words, topologically, a sequence is a connected oriented $1$-manifold with boundary. The natural next step are $2$-simplices, which matches the representation of ternary structures by $3$-graphs\footnote{Note that an oriented $3$-hyperedge with boundary is a $2$-simplex.} used in previous sections. We could define ternary sequences as simplicial complexes triangulating $2$-manifolds with the same general topological properties as ordinary sequences. Under this approach, the notion of ternary sequence is far from unique since there exist many classes of topologically inequivalent connected oriented 2-manifolds with boundary. A disk may serve as the simplest obvious template for the topology of ternary sequences, however other non-trivial topologies will certainly play a significant role in a complete theory of ternary sequentiality.

\section{Ternary Science} \label{ternaryscience}

In this section we present a broad spectrum of phenomena that display explicitly ternary or higher qualities. Our intention is merely to illustrate the far-reaching scope of the concept of arity across disciplines. Depending on the field of study, the higher arity perspective may provide suggestive new approaches, help re-conceptualize established knowledge or simply offer a neat recollection of anecdotal information.

\subsection{Biology} \label{bio}

\textbf{Ecology.} The need to account for higher order interactions between species has long been noted by ecologists \cite{abrams1983arguments,billick1994higher,wootton1994pieces,peacor2004dependent,mayfield2017higher,levine2017beyond,letten2019mechanistic}. More specifically, detailed accounts of 3-way symbiotic relations have been documented \cite{marquez2007virus,roopin2011amonia}.\newline

\textbf{Macromolecule Function.} The biological function of many molecular complexes can often be linked to low-arity interactions between structural sub-units \cite{phizicky1995protein,nooren2003diversity,bertoni2017modeling}. Two explicitly ternary examples of molecular complexes are porins, found in human cells \cite{song2015tnf}, and the triple-helix configuration of DNA strands, which is less stable than the double-helix but still occurs spontaneously in nature \cite{lestienne2011priming}.\newline

\textbf{Gene Expression.} Single discernible traits of a phenotype are often the result of complicated interactions of several gene expression processes \cite{bendor1999clustering,han2020heterogeneity}. Typical approaches to gene co-expression employ binary network methods \cite{zhao2010coexpression,pardo2020cogent}, however it has recently been noted that higher order gene interactions also play an important role \cite{taylor2015genetic}.\newline

\textbf{Organism Anatomy.} Higher order models of tissue connectivity have become increasingly popular in recent years \cite{klamt2009hypergraphs,bordbar2011multi}. The structure of slime mould outgrowths \cite{jabr2012brainless}, plant vascular systems and cognition \cite{calvo2020plants} or the function of glial cells in the brain \cite{allen2009glia} appear as promising applications of higher order tissue connectivity. Furthermore, in nervous systems, whose study has been dominated by binary network methods \cite{bullmore2011brain,van2016comparative}, there is mounting evidence for higher order connectivity and structure \cite{yu2011higher,petri2014homological,chambers2016synaptic,reimann2017cliques,tolokonnikov2019convolution,kilic2022cascades}.\newline

\textbf{Epidemiology.} The propagation of infectious diseases is typically modelled by diffusion mechanisms on binary networks \cite{lloyd2007models}, however new higher order approaches based on simplicial complexes \cite{matamalas2020phase} and hypergraphs \cite{bodo2016epidemic,higham2021epidemic} are becoming increasingly popular. 

\begin{figure}[h]
\centering
\includegraphics[scale=0.5]{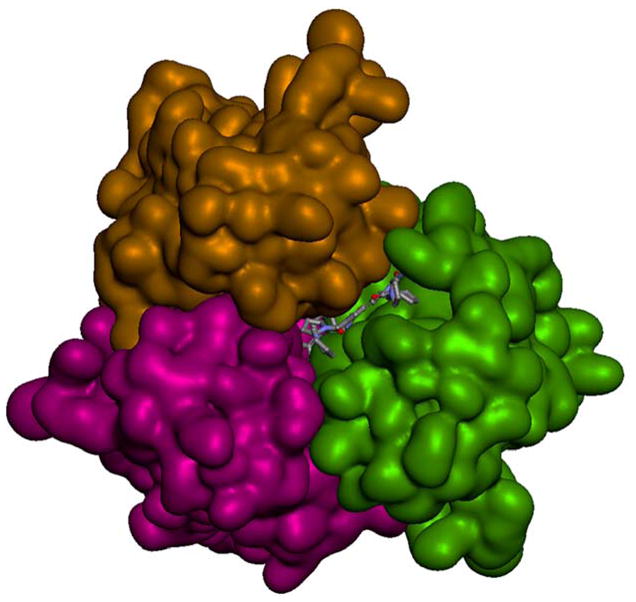}
\includegraphics[scale=0.27]{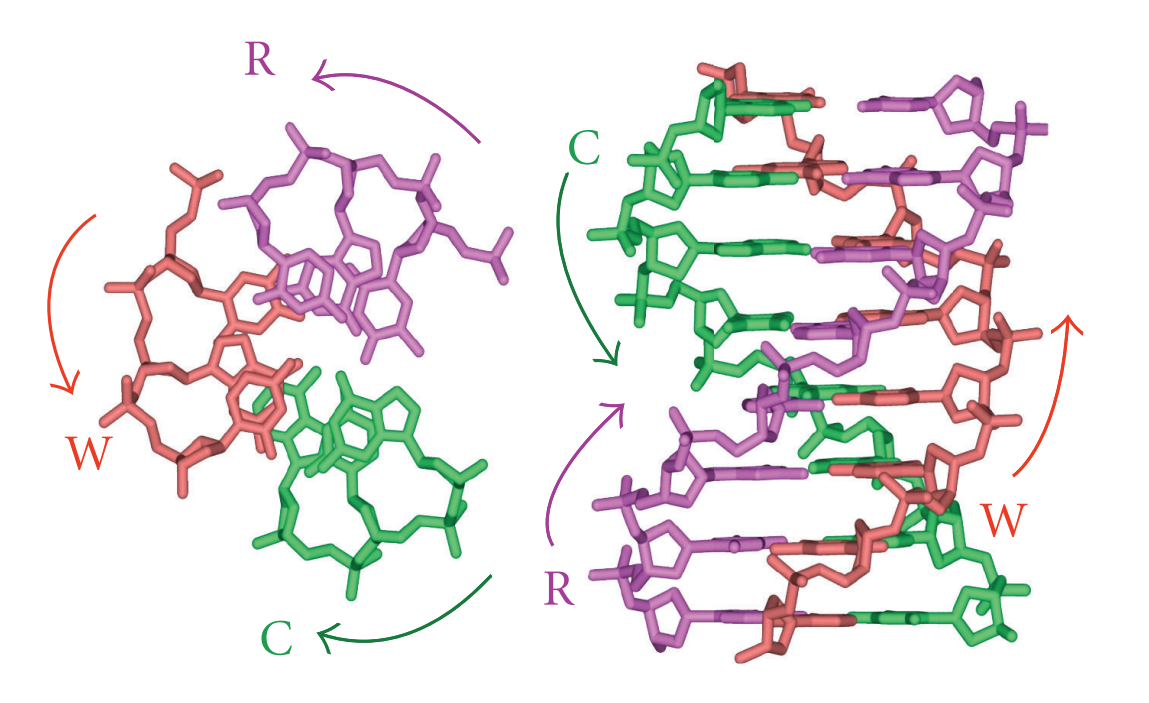}
\caption{\textbf{Biochemical Aritons}: (left) structure of the CD40L (CD154) protein trimer found in human cells \cite{song2015tnf}; (right) triple helix structure of R-form DNA strands \cite{lestienne2011priming}.}
\label{biochem}
\end{figure}

\subsection{Physics} \label{phys}

\textbf{Multi-Body Systems.} Beginning with the classical 3-body problem \cite{koon2000dynamical}, there is a long documented history of the difficulties faced in the study of multi-particle interactions and higher order effects \cite{doran1971multipole,cohen1993fifty,dobnikar2002many,battiston2021physics}. Here we list a few fascinating examples of ternariton interactions in classical, quantum, statistical and chemical systems: a ternary Boltzmann gas equation \cite{ampatzoglou2021rigorous}, three-body interactions in colloids \cite{dobnikar2004three,hynninen2004effect}, chemical compounds that display Borromean-like molecular structure \cite{cantrill2005borromean,meyer2010borromean,li2010borromean,baas2012chemical,pan2014borromean}, Borromean-like interaction of light nuclei \cite{zhukov1993borromean,bertulani2007borromean,betan2017cooper} and Efimov states of trios of particles \cite{cornelius1986efimov,fedorov1994efimov,hammer2010efimov}.\newline

\textbf{Nuclear Interaction.} The study of the strong nuclear force is noteworthy both in terms of the influence it has had in the development of many-particle formalisms \cite{kinoshita1954collective} and the marked ternary flavour of quark models \cite{gellmann1972quarks,fritzsch1973advantages,nambu1974three,hendry1978quarks}. Furthermore, the search for quantization schemes to derive such models led to the first use of an explicitly ternary algebra in theoretical physics \cite{nambu1973generalized,takhtajan1994generalized}.\newline

\textbf{Mathematical Physics.} This discipline deserves our special attention since most of our mathematical discussion on ternary structures in \ref{ternarymaths} is directly motivated by open problems in theoretical and mathematical physics. Higher arity ideas appear in mathematical physics in four different flavours:
\begin{itemize}
    \item \textbf{n-ary Algebras.} Since the first uses of cubic and ternary structures in the mid 20$^{\text{th}}$ century, particularly in nuclear physics as discussed above \cite{nambu1973generalized,nambu1974three}, n-ary algebras have appeared in physics-inspired mathematical research as n-ary associative algebras \cite{michor1998n,kerner2008ternary}, n-ary Lie algebras (also known as Filippov algebras) \cite{filippov1985n,bagger2008gauge,azcarraga2010nary,azcarraga2011contractions}, cubic matrix algebras \cite{kerner1997cubic,abramov2009algebras}, $\mathbb{Z}_3$-graded algebras \cite{celik2016some,abramov2017algebras}, ternary Clifford algebras \cite{abramov1996z_3,cerejeiras2021ternary}, and Nambu-Poisson structures \cite{nakanishi1998nambu,ibanez1999leibniz,grabowski2000filippov,debellis2010quantized}.
    \item \textbf{Higher Field Theories.} The application of modern category theory to theoretical physics has led to the development of higher gauge theories \cite{baez2011invitation}, higher algebras in supergravity \cite{ravera2021hidden}, homotopical extensions of quantum field theories \cite{schreiber2021higher}, homotopical pregeometric theories   \cite{arsiwalla2021pregeometric,arsiwalla2021homotopies}, higher topos theory \cite{lurie2009higher}, polycategory-based theories \cite{blanco2020bifibrations} and operad-based theories \cite{benini2021operads}. Although, as we pointed out in  \ref{cat}, most of these higher categorical notions are still firmly based on binary morphisms, they represent a general trend towards higher order concepts.
    \item \textbf{Many-Particle Systems.} The mathematical modelling of classical many-particle systems has led to the use of higher arity devices such as higher generalization of the Boltzmann equation \cite{ampatzoglou2021rigorous} or multi-particle instantaneous interactions \cite{dobnikar2002many}. On the quantum side, the most prominent appearance of higher arity ideas is in the phenomenon of higher order entanglement \cite{zeilinger1992higher,you2020higher}.
    \item \textbf{Generalizations of Duality.} The notion of duality appears both as a fundamental aspect of the mathematical underpinnings of physics, e.g. dual vector spaces or duality functors, as well as a tool to capture equivalences between two seemingly different physical theories, e.g. dualities in string theory \cite{alvarez1994some}. Recent research has found ternary analogues of these two notions of duality: a triality of vector spaces, generalizing the notion of dual pairing, in the context of the study of the $SO(8)$ symmetry group \cite{sridharan1994compositions,lounesto2001octonions,sapio2001spin,borghese2013triality,} and a triality of equivalent theories in physics, algebraic geometry and algebraic topology \cite{sati2021mysterious}.
\end{itemize}

\begin{figure}[h]
\centering
\includegraphics[scale=0.8]{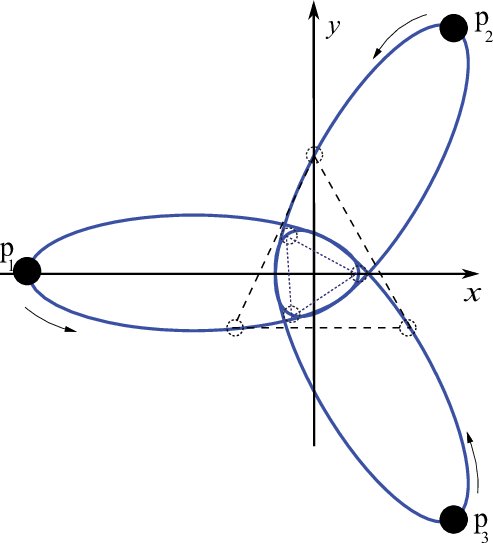}
\includegraphics[scale=0.32]{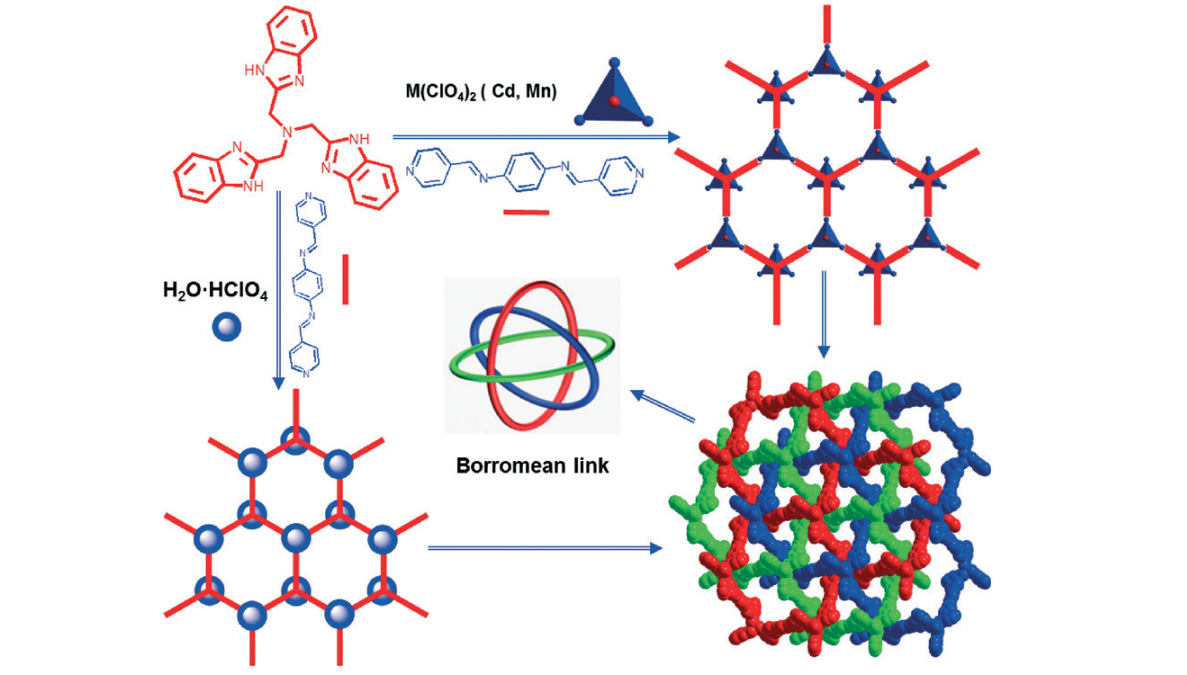}
\caption{\textbf{Physical Aritons}: (left) a planar solution of the gravitational 3-body problem, first discovered by J.L. Lagrange \cite{lagrange1772essai}, where three particles of identical mass always remain in the vertices of an equilateral triangle as they orbit the centre of mass (image credit: \cite{rivera2013periodic}); (right) hydrogen-bonding hcb-networks forming a 3-fold interlocking network of Borromean links \cite{pan2014borromean}.}
\label{physical}
\end{figure}

\subsection{Cognition} \label{cog}

\textbf{Colour Perception.} The particulars of human colour perception \cite{boynton1990human,king2005human} offer an intuitive source of ternary phenomena in the form of the experience of whiteness as emergent from the simultaneous combination of red, green and blue.\newline

\textbf{Flavour Perception.} In the subjective experience of taste, it is common to abstract and name a flavour from a complicated mix of chemical components, each having their distinctive individual flavour \cite{prescott2004psychological}. Interestingly, recent research on the ion channels involved in the nervous impulse for flavour and odour perception has shown that this higher arity is present even at the structural molecular level: several molecules can engage with  multiple binding sites of a single ion channel simultaneously \cite{del2021structural}.\newline

\begin{figure}[h]
\centering
\includegraphics[scale=0.57]{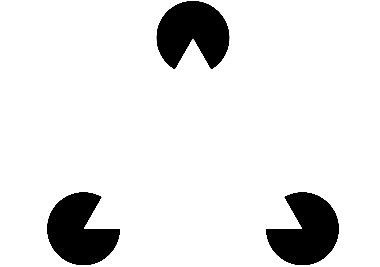}
\includegraphics[scale=0.3]{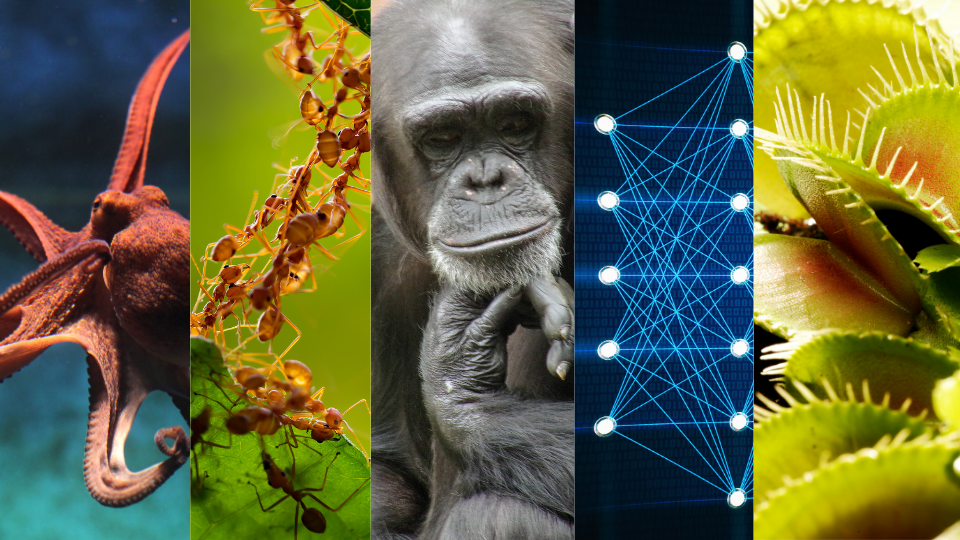}
\caption{\textbf{Cognitive Aritons}: (left) the illusion of a triangle resulting from the organization of three  partially cut-out disks, first illustrated by psychologist G. Kanizsa \cite{macpherson2017kanizsa}; (right) several systems of varying material substrate involving many interacting parts seem to give rise to emergent cognitive and/or conscious phenomena (image credit: Canva picture stock).}
\label{gestalt}
\end{figure}

\textbf{Gestalt Psychology.}  The birth of gestalt psychology in the early 20th century was an attempt to describe certain phenomenological aspects of visual perception such as `gist' or `context' that could not be adequately explained by merely taking into account properties of individual components of a stimulus \cite{wertheimer1938gestalt,koffka1935principles}. This led gestalt theorists to operationalize the idea of ``the whole being greater than the sum of its parts'' \cite{smith1988foundations}. This notion of holism proposed in gestalt psychology precisely alludes to higher structures and organizational principles in perception. From this perspective, illusions such as the Kanizsa triangle (Figure \ref{gestalt}) are interpreted as perceptually emergent phenomena resulting from the global organization of the parts comprising the stimulus. It has been reported that the specific processes that explain the emergence of these phenomena range from local feature detection to global strategies of perceptual organisation  \cite{spillmann1995phenomena}. A complete description of such perceptual phenomena will arguably require notions of higher arity.  \newline

\textbf{Conceptual Spaces.} A conceptual space is a non-symbolic cognitive model for higher-level associative representation of mental constructs and concept learning \cite{gardenfors2004conceptual}. Conceptual spaces consist of a number of quality dimensions derived from perceptual mechanisms and based on geometric similarity measures. For example, the concept of taste is represented as a geometric shape, a tetrahedron, whose vertices correspond to the bases tastes (sweet, sour, bitter, saline -- associated to 4 types of receptors). Then any other of a plethora of representable tastes are points or regions on this geometric structure. More recently, this has been extended to a cognitive theory of semantic representations  \cite{gardenfors2014geometry}. This framework heavily relies on higher order structures and would immensely benefit from a mathematical theory of higher arity for classifying mental representations.\newline

\textbf{Consciousness.} Considerations of holism similar to gestalt theory have also influenced the neuroscientific study of consciousness, where a particular informational complexity of the brain, known as integrated information, has been postulated as a measure of consciousness \cite{tononi1998consciousness,tononi2004information}. According to the Integrated Information Theory of consciousness, conscious experience is quantified in terms of the integrated information generated by the brain as a whole over and above the information generated individually by any of its parts \cite{oizumi2014phenomenology,arsiwalla2016high,arsiwalla2016global}. This integrated information captures the processing complexity associated to simultaneous integration and differentiation of the brain's structural and dynamical motifs at all architectural scales   \cite{koch2016neural,arsiwalla2017morphospace,arsiwalla2018measuring}. The qualia of consciousness are thus represented as informationally irreducible structures which strongly allude to notions of higher aritons \cite{tononi2016integrated}.

\begin{figure}[h]
\centering
\includegraphics[scale=0.52]{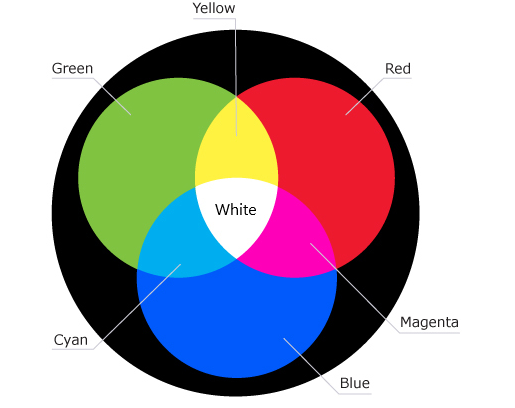}
\includegraphics[scale=0.29]{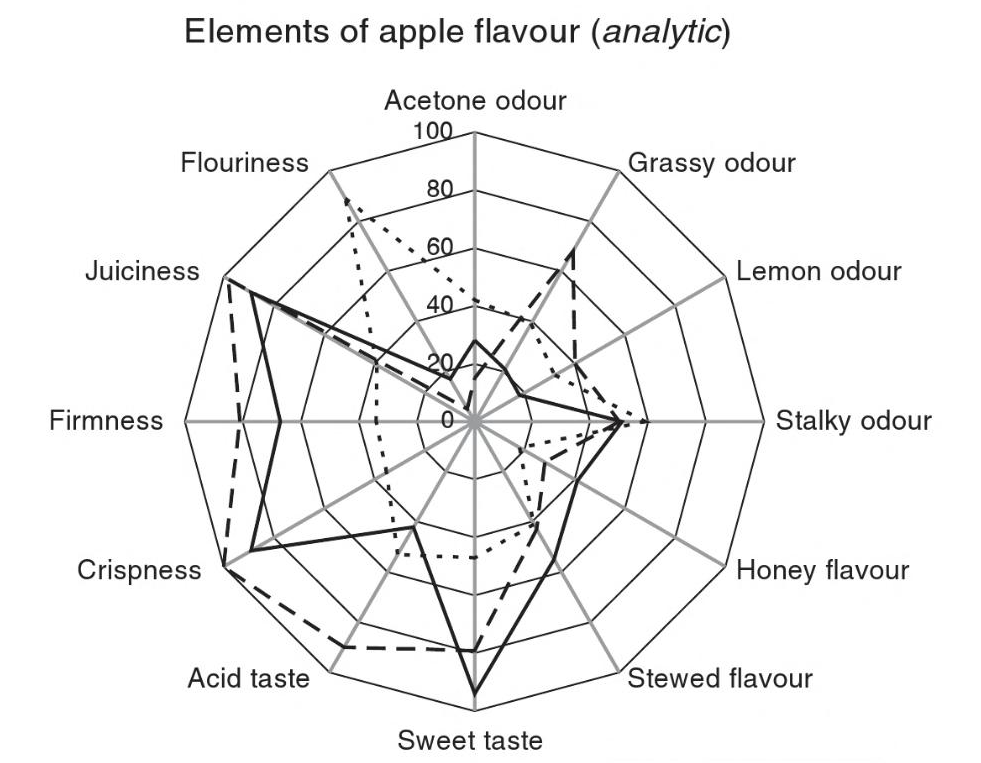}
\caption{\textbf{Perceptual Aritons}: (left) colours perceived by humans are aritons of primary colours, particularly, white is a ternariton under the additive colour mixing operation; (right) foods are usually conceptualized as a single flavour, however their components have distinct flavours that can be conceptualized individually \cite{prescott2004psychological} making food flavours into high-arity relations.}
\label{perceptual}
\end{figure}

\subsection{Computation} \label{comp}

\textbf{Hypernetworks.} Data scientists have long been noting the shortcomings of binary networks. Abundant research has been published advocating the use of hypergraphs and higher-order networks in data science and computation \cite{benson2016higher}, \cite{battiston2020networks}. The need for higher-order networks and models have also been proposed in the artificial neural network and machine learning communities \cite{zhou2006learning}. More recently,  hypernetwork architectures have been used for addressing the weight-sharing problem across layers of deep convolutional and recurrent networks \cite{ha2016hypernetworks}.   \newline

\textbf{Generalized Logic.} The core limitations of ordinary formal languages when it comes to describing higher arity arise from (i) the sequential nature of most symbolic representations (typically, strings of characters) and (ii) from the use of binary truth values. Progress concerning the first point has been made via the use of symbolic reasoning systems that exploit 2- and 3-dimensional diagrams to articulate formal reasoning. In particular, automated reasoning systems using string diagram rewriting methods make use of hypergraph structures formalized as tensor network diagrams \cite{bonchi2020string}. As for the second limitation, potential higher arity generalizations have been formulated as formal systems equipped with multiple symmetrical truth values. Specific examples of this case are type theories and their associated semantics that do away with the law of the excluded middle  \cite{program2013homotopy}.  \newline

\textbf{Distributed Processing.} Modern computer science is undergoing new paradigm shifts from single input-output-based sequential computation to distributed computation \cite{attiya2004distributed} on one hand, and  multicomputation-based processes \cite{wolfram2002new,wolfram2020class,wolfram2021multicomputation}, on the other. While the former is concerned with distributed algorithms involving simultaneous message-passing, \cite{lynch1996distributed}; the latter is concerned with all possible threads of computation in parallel, requiring a model of an observer to determine states \cite{wolfram2021multicomputation}. Both of these paradigms do away with sequential processes and offer new testing grounds for higher arity representation of computational processes. Recent examples utilizing the above paradigms include distributed consensus networks and blockchain systems \cite{wolframdconsensus}.  \newline

\textbf{Diagrammatic Quantum Computing.} In recent years, Categorical Quantum Mechanics arose as a diagrammatic formalism for describing quantum processes and protocols \cite{abramsky2009categorical}. This uses 2-dimensional category-theoretic string diagrams -- as opposed to the conventional sequential representation of operator actions -- for formal reasoning within systems of parallel quantum processes \cite{coecke2018picturing}. A particular class of this diagrammatic quantum process algebra, called ZX-Calculus, is used to represent and reason within quantum circuits \cite{coecke2011interacting}. This diagrammatic quantum formalism in fact gives a representation of quantum tensor networks and has been extremely successful at addressing problems related to quantum circuit simplification and optimization \cite{backens2014zx,jeandel2018complete,duncan2020graph,kissinger2019pyzx,gorard2021fast}. This framework of quantum computing naturally allows for higher arity operators and has recently been related to a class of multiway rewriting systems \cite{gorard2020zx,gorard2021zx}.

\subsection{Economics} \label{eco}

\textbf{Multiplayer Games.} The competitive paradigm of dyadic games has been extensively studied \cite{myerson1997game,tadelis2013game}, however, the mere addition of a third competitor introduces diplomacy and cooperation effects which dramatically increase the complexity of the strategy space that the players navigate \cite{shoham2009multiagent,camerer2011behavioral}. A multiplayer game that requires inputs from all players for its progression can be naturally regarded as a higher arity relation between agents. Although relatively understudied when compared to 2-player game theory, multiplayer game theory is already finding some promising applications in evolutionary theory  \cite{gokhale2014evolutionary}, linguistics \cite{havrylov2017emergence} and machine learning \cite{albrecht2018autonomous}.\newline

\textbf{Concurrent Competition.} The dynamics of competition, particularly in economic and financial systems, have been typically modelled successfully with binary networks \cite{markovich2009winning,braha2011corporate}, though in recent years there has been an increasing number of studies that approach the phenomenon of competition with higher order methods \cite{amato2017opinion,letten2019mechanistic}.\newline

\textbf{Tokenomics.} Digital and cryptographic currency technologies \cite{nakamoto2008bitcoin,mukhopadhyay2016brief} offer a flexible template where generalized paradigms of transaction, e.g. requiring the concurrent participation of three parties in any exchange, can be easily implemented, deployed and studied.

\begin{figure}[h]
\centering
\includegraphics[scale=0.25]{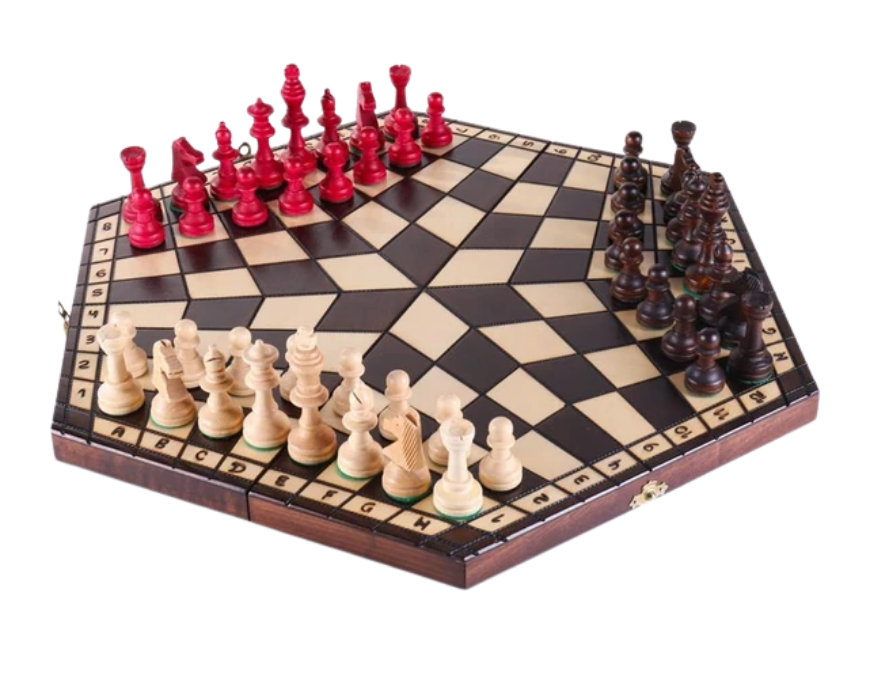}
\includegraphics[scale=0.2]{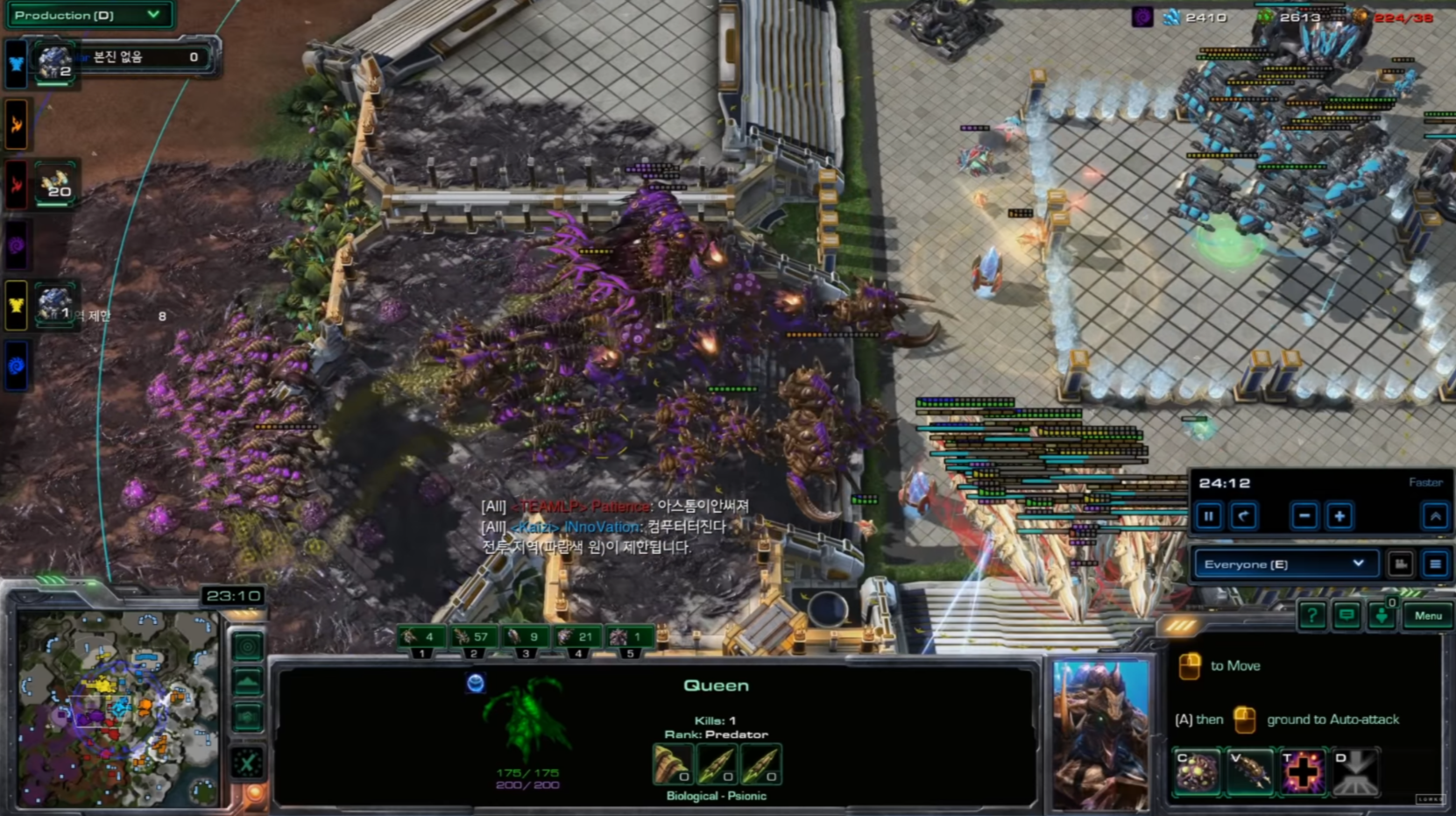}
\caption{\textbf{Competitive Aritons}: (left) 3-player chess, an example of a ternary turn-based strategy game with perfect information; (right) 3-player free-for-all StarCraft 2 map, an example of a ternary real-time strategy game with imperfect information. These examples suggest an extension of the research that led to the groundbreaking AIs AlphaZero \cite{silver2018general} and AlphaStar \cite{arulkumaran2019alphastar} by simply including 3-player games at the stage of agent training.}
\label{games}
\end{figure}

\subsection{Sociology} \label{soc}

\textbf{Collaboration.} Human societies are known to be a rich source of instances of higher order organization \cite{benson2016higher}. Collaborative activities, such as scientific research \cite{newman2001structure} or governance \cite{torfing2012governance}, give rise to natural aritons of concurrent human participation.\newline

\textbf{Association.} At a smaller scale, relationships between individuals offer significant instances of low arity in the form of affectional bonds \cite{ainsworth2006attachments} and small group dynamics \cite{homans2017human}. Romantic and sexual associations are perhaps particularly interesting given the observed frequency disparity between pairs and multi-party arrangements \cite{balzarini2019demographic}, which we may hypothesize as a socio-sexual analogue of the predominance of binary structures in intellectual disciplines.\newline

\textbf{Identity.} Belonging and group recognition in complex societies often result in the intersection of multiple notions of identity on a single individual \cite{simon2008identity,purdie2008intersectional}. Higher arity concepts might become helpful when navigating complex social environments by reducing the cognitive load on the individual when dealing with conflict.

\subsection{Arts} \label{art}

\textbf{Music.} Western classical music composition \cite{scagliarini2021quantifying} and repertoire \cite{park2015topology} has been recently investigated under the lens of higher order information. From a practitioner's point of view, harmonic relations between tones in a tuning system are natural examples of higher arity \cite{tymoczko2006geometry,hofmann2008virtual,chan2019science}. In this picture, chords can be regarded as harmony aritons and the harmonic development of a piece of music can be accurately represented as a path in the harmony hypergraph.\newline

\textbf{Architecture.} The rigidity properties of triangles makes them commonplace as structural elements in construction and architecture \cite{ross1994triangles}. From the point of view of equilibrium mechanics, ternary structures appear as the manifestation of a basic and intuitive statics fact: a tripod is the only piece of furniture that never wobbles.\newline

\textbf{Ritual and Spirituality.} Multiple human traditions have captured cultural and mystical beliefs in ternary concepts such as trios of entities or objects with some form of ternary symmetry \cite{pogliani2005intriguing}. Two prominent examples are the concept of the christian holy trinity \cite{letham2004holy} and the beliefs attached to symbols such as the triquetra or the valknut in pre-christian Norse and Germanic culture \cite{andren2011old}.

\begin{figure}[h]
\centering
\includegraphics[scale=0.135]{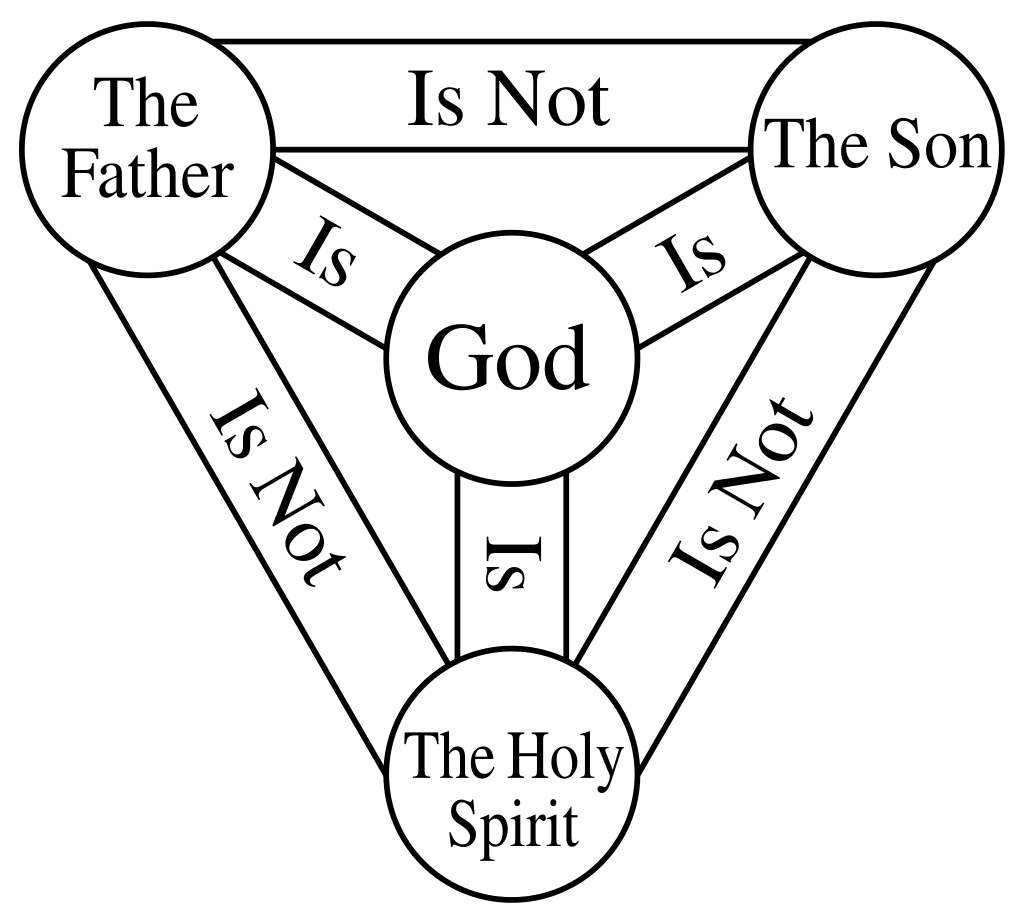}
\includegraphics[scale=0.322]{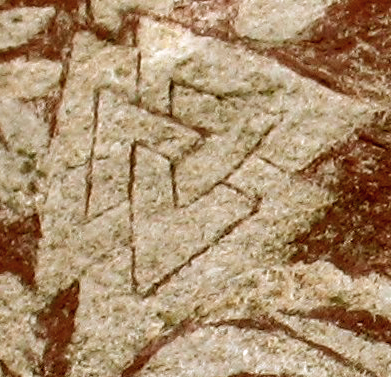}
\includegraphics[scale=0.36]{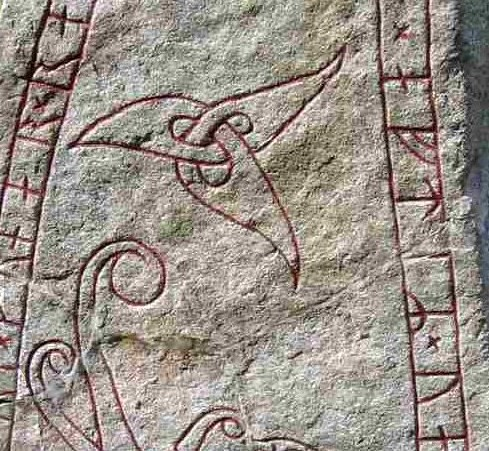}
\caption{\textbf{Mystical Aritons}: (left) a diagram representing the non-transitivity properties of the christian holy trinity; (centre) a valknut symbol on the Stora Hammars stone found in Gotland island, Sweden, dated around the year 600, this instance of the valknut displays the same topology as the Borromean rings; (right) a triquetra symbol on the U937 Funbo runestone found in the Uppland province, Sweden, dated around the year 1100, the triquetra displays the topology of a trefoil knot (image credits: Wikipedia).}
\label{mystical}
\end{figure}

\section{Towards Transdisciplinary Higher Arity Science} \label{concl}

Although higher arity methods are becoming increasingly popular across disciplines and higher order phenomena are something of a hot topic among some research communities these days, most approaches are often limited to theories and methodologies that merely extend conventional binary ideas. We believe that higher arity concepts and techniques should be embraced in their own right, challenges notwithstanding, as they present immense potential for scientific discovery and mathematical creativity. Indeed, the many natural phenomena and mathematical structures we have discussed in this article offer ample evidence for the existence of a rich universe of genuinely higher arity ideas waiting to be explored.\newline

In the early stages of any field of enquiry, particularly those that push the boundaries of human intuition as higher arity does, it is paramount that theoretical accounts do not venture too far away from the phenomena they attempt to describe and that preliminary models remain in a tight feedback loop with experiments. Successful development of higher arity science in the intellectual landscape of the 21$^{\text{st}}$ century would require the merging of formal and empirical disciplines and the blurring of boundaries between mathematics, science and computation. The transdisciplinary nature of higher arity ideas is perhaps best illustrated by the potential application of irreducible arity as a measure of atomic complexity or arity cohomology as a measure of emergent behaviour.\newline

The prospect of these future lines of research becomes all the more exciting when we consider how poorly we currently understand the simplest forms of higher arity. Beyond the well-documented history of the difficulties faced when modelling systems with ternary or higher interactions, even the most elementary examples of ternary structures behave in unexpected ways that cannot be easily grasped as direct generalizations of familiar binary phenomena. This intangible conceptual barrier limiting all kinds of intellectual disciplines is what we may call an \emph{unthought frontier of science}. The present article should be taken as a declaration of our determination to venture into uncharted territories and a rallying call for others to join what promises to be a thrilling expedition.

\section*{Author Contributions}

Carlos Zapata-Carratal\'a is the lead author of this article. Carlos is currently coordinating a small interdisciplinary collaboration team to pursue research on higher arity science as outlined in this article. Xerxes D. Arsiwalla is the second author and has contributed substantially to the research, development and writing of Sections \ref{cubix}, \ref{rel}, \ref{cat},  \ref{phys}, \ref{cog},   \ref{comp} as well as overall revisions throughout the paper.

\section*{Acknowledgements}

We would like to thank Irida Altman, Nils Baas, Tali Beynon, Jos\'e Figueroa-O'Farrill, Prathyush Pramod, Pau Enrique, \'Alvaro Moreno, Rafael N\'u\~nez, Stephen Wolfram, Mundy Otto Reimer and Silvia Butti for all the feedback and eye-opening conversations.

\printbibliography

\end{document}